\documentclass[11pt]{article}
\usepackage{amscd}
\usepackage{verbatim}
\usepackage{amssymb}
\usepackage{amsmath}
\usepackage{amsthm}
\usepackage{finalT}

\newcommand{\id}{\mathord{{\mathrm 1}\kern-0.27em{\mathrm I}}\kern0.35em}

\newcommand{\Half}{\ensuremath{\textstyle\frac{1}{2}}}
\newcommand{\Quarter}{\ensuremath{\textstyle\frac{1}{4}}}

\newcommand{\enorm}[1]{|#1|}
\newcommand{\norm}[1]{\|#1\|}

\newcommand{\ipe}[2]{ ( #1 | #2 )}

\newcommand{\ip}[2]{ \langle #1 | #2 \rangle}

\newcommand{\At}{\tilde{A}}
\newcommand{\Dc}{\mathcal{D}}

\newcommand{\gt}{\tilde{g}}
\newcommand{\htl}{\tilde{h}}

\newcommand{\Ht}{\tilde{H}}
\newcommand{\Nbb}{\mathbb{N}}
\newcommand{\Rbb}{\mathbb{R}}
\newcommand{\ttl}{\tilde{t}}

\newcommand{\Co}{\text{C}_{0}^{\infty}}

\newcommand{\Je}{J_{\epsilon}}
\newcommand{\ue}{u_{\epsilon}}
\newcommand{\ep}{\epsilon}
\newcommand{\Jde}{J^{\dagger}_{\epsilon}}
\newcommand{\ute}{\tilde{u}_{\epsilon}}

\begin{document}
\vskip 0 true cm \flushbottom
\begin{center}
\vspace{24pt} { \large \bf Asymptotically Flat Ricci Flows} \\
\vspace{30pt}
{\bf T.A. Oliynyk}$^{\dag}$ \footnote{todd.oliynyk@aie.mpg.de},
{\bf E Woolgar}$^{\ddag}$ \footnote{ewoolgar@math.ualberta.ca} %%

\vspace{24pt} %%
{\footnotesize $^\dag$ Max-Planck-Institut f\"ur
Gravitationsphysik (Albert Einstein Institute), Am M\"uhlenberg 1,
D-14476 Potsdam, Germany.\\
$^\ddag$ Dept of Mathematical and Statistical Sciences,
University of Alberta,\\
Edmonton, AB, Canada T6G 2G1. }
\end{center}
\date{\today}
\bigskip

%\smallskip
%\begin{center}\textbf{NOTE}: This paper has been accepted for publication
%in Comm. Anal. Geom. \end{center}
%\smallskip

\begin{center}
{\bf Abstract}
\end{center}
\noindent We study Ricci flows on ${\mathbb R}^n$, $n\ge 3$, that
evolve from asymptotically flat initial data. Under mild conditions
on the initial data, we show that the flow exists and remains
asymptotically flat for an interval of time. The mass is constant in
time along the flow. We then specialize to the case of rotationally
symmetric, asymptotically flat initial data containing no embedded
minimal hyperspheres. We show that in this case the flow is
immortal, remains asymptotically flat, never develops a minimal
hypersphere, and converges to flat Euclidean space as the time
diverges to infinity. We discuss the behaviour of quasi-local mass
under the flow, and relate this to a conjecture in string theory.

%\maketitle

\setcounter{equation}{0}
\newpage

\section{Introduction}
\setcounter{equation}{0}

\noindent The Ricci flow
\begin{equation}
\frac{\partial g_{ij}}{\partial t} = -2 R_{ij}\ . \label{eq1.1}
\end{equation}
was first introduced in the mathematics literature by Richard
Hamilton \cite{Hamilton1} in 1982. Almost immediately, it was
applied to the classification problem for closed 3-manifolds and
much subsequent work in the subject in the intervening 25 years has
been focused on this application, culminating in the recent
celebrated results of Perelman \cite{Perelman}.

By contrast, Ricci flow on noncompact manifolds has received
somewhat less attention. Of course, structures on noncompact
manifolds, such as Ricci solitons, are relevant to the compact case,
and this has been to now an important motivation for work on the
noncompact case. The case of asymptotically flat Ricci flow has
remained virtually untouched (nontrivial solitons do not occur in
this case \cite{OSW1}).

But physics provides considerable motivation to study the
asymptotically flat case. Our interest in it arises out of a
conjectural scenario in string theory. Equation (\ref{eq1.1}) is the
leading-order {\it renormalization group flow equation} for a
nonlinear sigma model that describes quantum strings propagating in
a background spacetime \cite{Friedan}.\footnote
{We ignore the dilaton since it can be decoupled from the metric in
renormalization group  flow.}
What is important to understand from this statement is that fixed
points of this equation provide geometric backgrounds in which the
low energy excitations of quantum strings can propagate (in the
approximation that radii of curvature are large and excitation
energies small relative to the so-called string scale).

The variable $t$ in renormalization group flow is not time: it is (a
constant times) the logarithm of the so-called renormalization
scale. However, there are conjectured relationships between
renormalization group flow and temporal evolution. A specific case
concerns tachyon condensation, the scenario wherein an unstable
string system is balanced at the top of a hill of potential energy
(for a review of tachyon condensation, see \cite{HMT}). The system
falls off the hill, radiating away energy in gravitational waves.
The system comes to rest in a valley representing a stable minimum
of potential energy. In open string theory, a more elaborate version
of this scenario involving the evaporation of a brane and the
formation of closed strings is now well understood, even
quantitatively. In closed string theory, much less is known but,
conjecturally, the fixed points of the renormalization group flow
equation (\ref{eq1.1}) are the possible endpoints of this evolution.
Sometimes it is further conjectured that time evolution in closed
string theory near the fixed points is determined by renormalization
group flow, and then $t$ in (\ref{eq1.1}) does acquire an
interpretation as a time.

Comparing both sides of this picture, we see that the radiation of
positive energy in the form of gravitational waves as the system
comes to rest in the valley should produce a corresponding decrease
in the mass of the manifold under the Ricci flow. This suggests that
we should endeavor to formulate and test a conjecture that mass
decreases under Ricci flow, at least if the initial mass is
positive.

The asymptotically flat case has a well-defined notion of mass, the
ADM mass, so this seems an appropriate setting in which to formulate
the conjecture. However, the metric entering the renormalization
group flow or Ricci flow in this scenario is not the full spacetime
metric, for which (\ref{eq1.1}) would not be even quasi-parabolic,
but rather the induced Riemannian metric on a suitable spacelike
submanifold \cite{GHMS}. Now ADM mass is conserved (between Cauchy
surfaces, and in the closed string scenario of \cite{GHMS}), even in
the presence of localized sources of radiation. This, we will see,
is reflected in the Ricci flow. The mass of $g$ will not change
during evolution by (\ref{eq1.1}). But if energy loss through
gravitational radiation occurs, then the {\it quasi-local mass}
contained within a compact region {\it should} change along the flow
to reflect this.\footnote
{We prefer not to discuss in terms of the Bondi mass, which would
require us to pass back to the Lorentzian setting which is not our
focus in this article. See \cite{GHMS} for a discussion in terms of
Bondi mass.}

In this paper, we focus first on the asymptotically flat case of
Ricci flow in general. Section 2 describes asymptotically flat
manifolds, with no assumption of rotational symmetry. Continuing
with the general asymptotically flat case, in Subsection 3.1 we
state and prove our short-term existence result Theorem 3.1, showing
that a general asymptotically flat data set on ${\mathbb R}^n$ will
always evolve under Ricci flow, remaining smooth and asymptotically
flat on a maximal time interval $[0,T_M)$. We will show that the ADM
mass remains constant during this interval, at least for
non-negative scalar curvature (i.e., the positive mass case, the
usual case of physical interest). Moreover, if $T_M < \infty$ then
the norm of the Riemann curvature must become unbounded as $t
\nearrow T_m$, just as in the compact case. We show this in
Subsection 3.2.

The short-term existence proof in Section 3 depends on detail
provided in the appendices. In Appendix A, we derive weighted
versions of standard Sobolev estimates such the Sobolev inequalities
and Moser estimates. We then use these estimates in Appendix B to
prove local existence and uniqueness in weighted Sobolev spaces for
uniformly parabolic systems.

We specialize to rotational symmetry in Section 4. In Section 4.1,
we pass to a coordinate system well suited to our subsequent
assumption that no minimal hyperspheres are present initially. We
show in Section 4.3 that this coordinate system remains well-defined
on the interval $[0,T_m)$. This is essentially a consequence of the
result, proved in Section 4.2, that no minimal hyperspheres develop
during the flow.

The absence of minimal spheres allows us to analyse the problem in
terms of a single PDE, the master equation (\ref{eq4.18}). From this
equation, we derive a number maximum principles that yield uniform
bounds on the curvature which allow us to conclude that
$T_M=\infty$. We obtain these principles in the first two
subsections of Section 5. Even better, we obtain not just uniform
bounds but decay estimates, from which we can prove convergence to
flat Euclidean space. Now given our assumptions, this is the only
Ricci-flat fixed point available. That is, the string theory
discussion above would lead one to conjecture that:

\bigskip
\noindent {\it When no minimal hypersphere is present, rotationally
symmetric, asymptotically flat Ricci flow is immortal and converges
to flat space as $t\to\infty$;}
\bigskip

\noindent and this is what we show. Though we have motivated this
conjecture from string theory for the case of positive initial mass,
we will prove that it holds whether or not the initial mass is
positive. This is our main theorem, proved in Subsection 5.3, which
states:

\begin{thm}\label{Thm1.1}
Let $\{x^i\}_{i=1}^n$ be a fixed Cartesian coordinate system on
${\mathbb R}^n$, $n\ge 3$. Let ${\hat g}=\hat{g}_{ij}dx^i dx^j$ be
an asymptotically flat, rotationally symmetric metric on ${\mathbb
R}^n$ of class $H^k_{\delta}$ with $k>n/2 +4$ and $\delta<0$. If
$(\Rbb^n,\hat{g})$ does not contain any minimal hyperspheres, then
there exists a solution $g(t,x) \in C^\infty( (0,\infty)\times
\Rbb^n)$ to Ricci flow (\ref{eq1.1}) such that
\begin{itemize}
\item[(i)] $g(0,x)=\hat{g}(x)$,
\item[(ii)] $g_{ij} - \delta_{ij} \in C^1([0,T],H^{k-2}_\delta)$
and $g_{ij}-\delta_{ij} \in C^1([T_1,T_2],H^\ell_\delta)$ for any
$0<T_1<T_2<\infty$, $0<T<\infty$, $\ell \geq 0$,
\item[(iii)] for each integer $\ell \geq 0$ there exists a constant
$C_\ell > 0$ such that
\eqn{curvest}{
\sup_{x\in\Rbb^n}|\nabla^\ell{\rm Rm}(t,x)|_{g(t,x)}\leq
\frac{C_\ell}{(1+t)t^{\ell/2}}
\quad \forall\; t>0\, ,
}
\item[(iv)]
the flow converges to $n$-dimensional Euclidean space ${\mathbb
E}^n$ in the pointed Cheeger-Gromov sense as $t\to\infty$, and
\item[(v)] if furthermore
$k>n/2 +6$, $\delta<\min\{4-n,1-n/2\}$, $\hat{R}\geq 0$, and
$\hat{R}\in L^1$, then the ADM mass of $g(t)$ is well defined and
$\text{\rm mass}(g(t)) = \text{\rm mass}(\hat{g})$ for all $t\geq
0$.
\end{itemize}
\end{thm}

When a minimal hypersphere {\it is} present initially, if the neck
is sufficiently pinched then we expect long-time existence to fail.
To see why, consider rotationally symmetric metrics on $S^n$. If
there is a sufficiently pinched minimal $(n-1)$-sphere, the
curvature blows up in finite time. This has been shown both
rigorously $(n\geq 3)$ \cite{AngKno} and numerically $(n=3)$
\cite{GI}. Our assumption of no minimal spheres in the initial data
is intended to prevent this. The ability to make this assumption and
to choose coordinates adapted to it is a distinct advantage of the
noncompact case. However, we also expect (based, e.g., on
(\cite{GI}) that for initial data with minimal hyperspheres that
have only a mild neck pinching, the flow will continue to exist
globally in time as well. Thus, when a minimal hypersphere is
present, we believe there would be considerable interest in
determining a precise criterion for global existence in terms of the
degree of neck pinching because of the possibility, raised in
\cite{GI}, that the critical case on the border between singularity
formation and immortality may exhibit universal features such as
those observed in critical collapse in general relativity
\cite{Choptuik}.

The constancy of the ADM mass in statement (v) is not at odds with
the conclusion that the flow converges to a flat and therefore
massless manifold. This constancy was also noted in \cite{DM} but we
draw different conclusions concerning the limit manifold, owing to
our use of the pointed Cheeger-Gromov sense of convergence of
Riemannian manifolds.\footnote
{The rotationally symmetric, expanding soliton of \cite{GHMS} can be
used to illustrate this phenomenon explicitly (albeit in 2
dimensions, whereas our results are for $n\ge 3$ dimensions). For
this soliton, one can easily compute the Brown-York quasi-local mass
on any ball whose proper radius is fixed in time and see that for
each such ball the quasi-local mass tends to zero as $t\to \infty$,
and the flow converges to Euclidean 2-space. But the mass at
infinity of the soliton (the deficit angle of the asymptotic cone in
2 dimensions) is a constant of the motion which can be set by
initial conditions to take any value.}
In Subsection 4.4 we define three different kinds of metric balls in
$({\mathbb R}^n,g(t))$, $n\ge 3$; balls of fixed radius, fixed
volume, and fixed surface area of the bounding hypersphere. To
clarify the behaviour of the mass in the limit $t\to\infty$, we
express the Brown-York quasi-local mass of these balls in terms of
sectional curvature and, by anticipating the decay rate for
sectional curvature derived in Section 5, show that these
quasi-local masses go to zero as $t\to\infty$, even though the ADM
mass, as measured at infinity, is constant. The picture is not
strongly dependent on the definition of quasi-local mass, of which
the Brown-York definition is but one among many. In rotational
symmetry in any dimension, the metric has only one ``degree of
freedom''. The study of the evolution of quasi-local mass then
reduces to the study of this single degree of freedom, no matter
which definition of quasi-local mass one prefers.\footnote
{The assumption of spherical symmetry in general relativity
precludes gravitational radiation, according to the Birkhoff
theorem. But on the string side of our scenario, the picture is one
of closed strings existing as perturbations that break the spherical
symmetry of the background metric (as well, we should include a
dilaton background field that modifies general relativity). Viewed
in the string picture, these perturbations create the radiation that
is detected as a change in the quasi-local mass of the spherically
symmetric Ricci flow.}

Although local existence, uniqueness, and a continuation principle
for Ricci flow on non-compact manifolds with bounded curvature are
known \cite{Shi,CZ}, it does not follow immediately from these
results that Ricci flow preserves the class of asymptotically flat
metrics. One of the main results of this paper is to show that Ricci
flow does in fact preserve the class of asymptotically flat metrics.
%%EW:Rewording so as not to imply that Dai-Ma and List precede us:
Independent of our work, Dai and Ma have recently announced that
they have also been able to establish this result \cite{DM}, as has
List in his recent thesis \cite{List}.
%%EW: Removed:
%%Recently this has been established independently in \cite{List,DM}.

Our approach to the problems of local existence, uniqueness,
continuation, and asymptotic preservation is to prove a local
existence and uniqueness theorem for quasi-linear parabolic
equations with initial data lying in a weighted Sobolev space, and
then use it to show that Ricci flow preserves the class of
asymptotically flat metrics.
%%EW: Removed words:
%%The main reason for taking
An important advantage of this approach rather than appealing to the
results of \cite{Shi,CZ,List,DM} is that we
%%want
obtain a local existence and uniqueness theorem on asymptotically
flat manifolds that is valid for other types of geometric flows to
which the results of \cite{Shi,CZ,List,DM} do not immediately apply,
% EW: Added:
and which are of interest in their own right.
For example, our local existence results contained in appendix B
combined with the DeTurck trick will yield local existence,
uniqueness, and a continuation principle for the following flows on
asymptotically flat manifolds: \eqn{SEF}{ \left.
\begin{array}{l}\partial_t g = - 2R_{ij} + 4 \nabla_i u \nabla_j u \\
\partial_t u = \Delta u
\end{array}
\right\} \quad \text{(static Einstein flow)}, \\
}
\eqn{RGF}{
\left.
\begin{array}{l}
\partial_tg_{ij} = - \alpha^{'}
\bigl(R_{ij} + \nabla_i \nabla_j\Psi +\Quarter
H_{jpq}H_{j}{}^{pq}\bigr) \\
\partial_t \Psi = \frac{\alpha^{'}}{2} (\Delta
\Psi - \vert \nabla \Psi \vert^2 +\vert H\vert^2) \\
\partial_t B_{ij}
= \frac{\alpha^{'}}{2}( \nabla^{k}H_{kij} - H_{kij}\nabla^{k}\Psi
) \quad ( H:= \text{d} B)
\end{array}
\right\} \quad \text{($1^{\text{st}}$ order sigma model RG flow),} }
and \eqn{RGF2}{
\partial_tg_{ij} = - \alpha^{'}
\bigl(R_{ij} + \frac{\alpha^{'}}{2}R_{iklm}R_j{}^{klm}\bigr) \quad
\text{($2^{\text{nd}}$ order sigma model RG flow with $B=\Phi=0$).}
} We note that the static Einstein flow has been previously
considered in the thesis \cite{List}. There a satisfactory local
existence theory on noncompact manifolds is developed and an also a
continuation principle for compact manifolds is proved.

The problem of global existence for rotationally symmetric metrics
on $\Rbb^3$ has previously been investigated in \cite{Ivey}. There
the assumptions on the initial metric are different than ours.
Namely, the initial metric in \cite{Ivey} has positive sectional
curvature and the manifold opens up as least as fast as a
paraboloid. Under these assumptions, it is shown that Ricci flow
exists for all future times and converges to either a flat metric
or a rotationally symmetric Ricci soliton.

Finally, throughout we fix the dimension of the manifold to be $n\ge
3$. As well, we usually work with the Hamilton-DeTurck form of the
Ricci flow
\begin{equation}
\frac{\partial g_{ij}}{\partial t} = -2R_{ij}+\nabla_i\xi_j
+\nabla_j\xi_i\ , \label{eq1.2}
\end{equation}
which is obtained from the form (\ref{eq1.1}) by allowing the
coordinate basis in which $g_{ij}$ is written to evolve by a
$t$-dependent diffeomorphism generated by the vector field $\xi$.

\bigskip
\noindent{\bf Acknowledgments.} We thank Suneeta Vardarajan for
discussions concerning the string theory motivation for this work.
EW also thanks Barton Zwiebach for his explanation of the rolling
tachyon. This work was begun during a visit by TO to the Dept of
Mathematical and Statistical Sciences of the University of Alberta,
which he thanks for hospitality. The work was partially supported by
a Discovery Grant from the Natural Sciences and Engineering Research
Council of Canada.

\section{Asymptotically flat manifolds}
\setcounter{equation}{0}

\noindent The definition of asymptotically flat manifolds that we
employ requires the use of weighted Sobolev spaces, which we will
now define. Let $V$ be a finite dimensional vector space with inner
product $\ipe{\cdot}{\cdot}$ and corresponding norm $\enorm{\cdot}$.
For $u\in L^p_{\text{loc}}(\Rbb^n,V)$, $1\leq p \leq \infty$, and
$\delta \in \Rbb$, the \emph{weighted $L^{p}$ norm} of $u$ is
defined by
\leqn{wLpdef}{ \norm{u}_{L^{p}_{\delta}} := \left\{\begin{array}{ll}
\norm{\sigma^{-\delta-n/p}\,u}_{L^p} &
\text{if $ 1\leq p < \infty$}\\
\\
\norm{\sigma^{-\delta}\,u}_{L^{\infty}} & \text{if $p=\infty$}
\end{array} \right. }
with
\begin{equation}
\sigma(x) := \sqrt{1+|x|^2}\ .
\end{equation}
The \emph{weighted Sobolev norms} are then given by
\leqn{wSobdef}{ \norm{u}_{W^{k,p}_{\delta}} := \left\{
\begin{array}{ll} \displaystyle{\Bigl(\sum_{|I|\leq k}
\norm{D^{I}u}^{p}_{L^{p}_{\delta-|I|} } \Bigr)^{1/p}} &
\text{if $1\leq p < \infty$} \\
\\ \displaystyle{\sum_{|I|\leq k}
\norm{D^{I}u}_{L^{\infty}_{\delta-|I|} }} & \text{if $p=\infty$}
\end{array}\right. }
where $k\in \Nbb_0$, $I = (I_{1},\ldots, I_{n}) \in \Nbb_{0}^{n}$ is
a multi-index and $D^{I} =
\partial_{1}^{I_{1}}\ldots\partial_{n}^{I_{n}}$. Here
$\partial_i = \frac{\partial\;}{\partial x^i}$ and
$(x^1,\ldots,x^n)$ are the standard Cartesian coordinates on
$\Rbb^n$. The weighted Sobolev spaces are then defined as
\eqn{wsobdef}{ W^{k,p}_{\delta} = \{\, u \in
W^{k,p}_{\text{loc}}(\Rbb^n,V) \, | \, \norm{u}_{W^{k,p}_{\delta}} <
\infty \,  \}\, .  }
Note that we have the inclusion
\leqn{inclusion}{ W^{k,p}_{\delta_1} \subset W^{\ell,p}_{\delta_2}
\quad \text{for $k\geq \ell$, $\delta_{1} \leq  \delta_2$}}
and that differentiation $\partial_i \:
:\:W^{k,p}_{\delta}\rightarrow W^{k-1,p}_{\delta-1}$ is continuous.
In the case $p=2$, we will use the alternative notation
$H^{k}_{\delta} = W^{k,2}_{\delta}$. The spaces $L^{2}_{\delta}$ and
$H^{k}_{\delta}$ are Hilbert spaces with inner products
\leqn{L2ip}{ \ip{u}{v}_{L^{2}_{\delta}} := \int_{\Rbb^{n}}
\ipe{u}{v} \sigma^{-2\delta-n}d^{n}x} and \leqn{Hkip}{
\ip{u}{v}_{H^{k}_{\delta}} := \sum_{|I|\leq k}
\ip{D^{I}u}{D^{I}v}_{L^{2}_{\delta-|I|} }\, , }
respectively.

As with the Sobolev spaces, we can define weighted version of the
bounded $C^k$ function spaces
$C^{k}_{b} := C^k(\Rbb^n,V)\cap W^{k,\infty}$ spaces. For a map
$u\in C^{0}(\Rbb^{n},V)$ and $\delta\in \Rbb$, let
\eqn{rsobdef3.1}{ \norm{u}_{C^{0}_{\delta}} := \sup_{x\in
\Rbb^{n}}|\sigma(x)^{-\delta}u(x)| \, . } Using this norm, we define
the $\norm{\cdot}_{C^{k}_{\delta}}$ norm in the usual way:
\eqn{rsobdef3.3}{ \norm{u}_{C^{k}_{\delta}} := \sum_{|I|\leq k}
\norm{\partial^{I}u}_{C^{0}_{\delta-|I|}} \, . } So then
\gath{rsobdef3.5}{
 C^{k}_{\delta} := \bigl\{\, u \in C^{k}(\Rbb^n,V)
\, | \, \norm{u}_{C^{k}_{\delta}} < \infty \: \bigr\} \, .
}

We are now ready to define asymptotically flat manifolds.
\begin{Def}
\label{asymdef} \mnote{[asymdef]} {\em Let $M$ be a smooth,
connected, $n$-dimensional manifold, $n\ge 3$, with a Riemannian
metric $g$ and let $E_R$ be the exterior region $\{\,\, x\in
\Rbb^n\,\,|\,\, |x|>R\}$. Then for $k>n/2$ and $\delta <0$, $(M,g)$
is {\em asymptotically flat of class} $H^k_\delta$ if
\begin{itemize}
\item[(i)] $g\in H^{k}_{\text{loc}}(M)$,
\item[(ii)] there exists a finite collection $\{U_{\alpha}\}_{\alpha=1}^{m}$
of open subsets of $M$  and diffeomorphisms $\Phi_\alpha : E_R
\rightarrow U_\alpha$ such that $M\setminus \cup_{\alpha}U_{\alpha}$
is compact, and
\item[(iii)] for each $\alpha\in \{1,\ldots,m\}$, there exists an $R>0$
such that $(\Phi_{\alpha}^{*}g)_{ij} - \delta_{ij} \in
H^k_\delta(E_R)$, where $(x^1,\ldots,x^n)$ are standard Cartesian
coordinates on $\Rbb^n$ and  $\Phi_\alpha^*g =
(\Phi_{\alpha}^{*}g)_{ij} dx^i dx^j$.
\end{itemize}
}
\end{Def}
The integer $m$ counts the number of asymptotically flat ``ends'' of
the manifold $M$. As discussed in the introduction, we are interested
in manifolds where $M\cong \Rbb^n$ and hence $m=1$. In this case, we can
assume that
$g=g_{ij}dx^i dx^j$ is a Riemannian metric
on $\Rbb^n$ such that
\leqn{falloff}{ g_{ij}-\delta_{ij}\, ,\; g^{ij}-\delta^{ij} \in
H^{k}_\delta }
where $g^{ij}$ are the components of the inverse metric, satisfying
$g^{ij}g_{jk} = \delta^{i}_k$. We note that results of this section
and Theorems \ref{LocA}, \ref{LocB}, and \ref{cont} of the next
section are are easily extended to the general case. We leave the
details to the interested reader.

In the following section, we will need to use diffeomorphisms
generated by the flows of time-dependent vector fields and also
their actions on the metric and other geometrical quantities.
Therefore, we need to understand the effect of composing a map in
$H^{k}_\delta(\Rbb^n,V)$ with a diffeomorphism on $\Rbb^n$.
Following Cantor \cite{Cantor}, we define
\eqn{diffdef}{ \Dc^{k}_{\delta} := \{\; \psi : \Rbb^n \rightarrow
\Rbb^n \, |\, \text{$\psi-\id \in H^{k}_\delta$, $\psi$ is
bijective, and $\psi^{-1}-\id \in H^{k}_\delta$} \} }
which is the group of diffeomorphisms that are asymptotic to the
identity at a rate fast enough so that the difference lies in
$H^k_\delta$. We will need to understand not only when composition
preserves the $H^k_\delta$ spaces but also when composition
$(\psi,u) \mapsto u\circ\psi$ is continuous as a map from
$\Dc^k_\delta \times H^k_\delta$ to $H^k_\delta$. In \cite{Cantor},
Cantor studied this problem under the assumption that $\delta \leq
-n/2$. He assumed $\delta \leq -n/2$ because that was what he needed
to prove the weighted multiplication lemma (see Lemma \ref{SobB}).
However, it is clear from his arguments that the proofs of his
results are valid whenever the multiplication lemma holds and
$H^{k}_\delta \subset C^{1}_b$. Therefore, by Lemmata \ref{SobA} and
\ref{SobB}, his results are valid for  $\delta \leq 0$.
\begin{thm}
{\emph{[Corollary 1.6,\cite{Cantor}]}} \label{Can3} \mnote{[Can3]}
For $k > n/2+1$ and $\delta \leq 0$, the map induced
by composition
\eqn{Can3.1}{
H^k_{\delta} \times \Dc^k_{\delta} \longrightarrow
H^k_{\delta} \: : \: (u,\psi) \longmapsto u \circ \psi
}
is continuous.
\end{thm}
Cantor also proved the following three useful results:
\begin{lem}
{\emph{[Lemma 1.7.2,\cite{Cantor}]}} \label{Can3B} \mnote{[Can3B]}
If $k > n/2+1$, $\delta \leq 0$, and $f$ is a $C^1_b$
diffeomorphism such that $f-\id \in H^{k}_\delta$
then $f \in \Dc^k_\delta$.
\end{lem}
\begin{thm}
{\emph{[Theorem 1.7,\cite{Cantor}]}} \label{Can3A} \mnote{[Can3A]}
For $k > n/2+1$ and $\delta \leq 0$, $\mathcal{D}^k_\delta$
is an open subset of
\eqn{Can3A.1}{
\{ \, f : \Rbb^n \rightarrow \Rbb^n\,|\, f-\id \in H^k_\delta \,\}\,.
}
\end{thm}
\begin{thm}
{\emph{[Theorem 1.9,\cite{Cantor}]}} \label{Can2} \mnote{[Can2]}
For $k > n/2+1$ and $\delta \leq 0$, $\Dc^k_\delta$ is a topological
group under composition and a smooth Hilbert manifold. Also, right
composition is smooth.
\end{thm}

The following proposition is a straightforward extension
of Cantor's work.
\begin{prop}
\label{Can4} \mnote{[Can4]}
If $k> n/2+1$, $\delta \leq 0$, and $u \in H^{k+\ell}_\delta$ $(\ell\geq 0)$
then
the map
\eqn{Can4A}{
\Dc^{k}_\delta \longrightarrow  H^k_\delta
\: : \: \psi \longmapsto u\circ \psi
}
is of class $C^\ell$.
\end{prop}
Using these results, it is not difficult to see that the proof of
Theorem 3.4 of \cite{EM} generalizes to the $H^k_\delta$ spaces with
the result being:
\begin{thm}
\label{EM1} \mnote{[EM1]} Suppose $\delta \leq 0$,
$k > n/2+2$ and
$X : (-\kappa,\kappa)\times \Rbb^n \rightarrow \Rbb^n$ $(\kappa >0)$
defines a continuous map
\eqn{EM1A}{X : (-\kappa,\kappa)\longrightarrow
H^{k+\ell}_\delta(\Rbb^n,\Rbb^n)\: : \: t \longmapsto
X(t,\cdot)\quad (\ell \geq 0)\, . }
Let $\psi_t$ denote the flow of the time dependent vector field
$X(t,x)$ on $\Rbb^n$ that satisfies $\psi_0 = \id$. Then there
exists a $\kappa_* \in (0,\kappa)$ such that $\psi_t$ $(\,t\in
(-\kappa_*,\kappa_*) \,)$ defines a $C^{1+\ell}$ curve in
$\Dc^k_\delta$.
\end{thm}

\section{Local Existence}
\setcounter{equation}{0}

\subsection{Existence of General Asymptotically Flat Ricci Flows}

\noindent We now prove a local existence result for Ricci flow on
asymptotically flat manifolds.

\begin{thm}
\label{LocA} \mnote{[LocA]}
Let $\hat{g}$ be an asymptotically flat metric of class $H^k_\delta$
with $\delta < 0$ and $k > n/2+3$. Then there exists a $T>0$ and
a family $\{g(t),t\in [0,T)\}$ of asymptotically flat metrics
of class $H^{k-2}_\delta$ such that
$g(0) = \hat{g}$,
\eqn{LocA1}{
g_{ij}-\delta_{ij}\, ,\; g^{ij}-\delta^{ij} \in
C^{1}([0,T),H^{k-2}_\delta)\, ,
}
and $\partial_t g_{ij} = -2R_{ij}$ for all $t\in [0,T)$. Moreover,
$g(t,x)\in C^\infty((0,T)\times M)$ and $g_{ij}-\delta_{ij}$,
$g^{ij}-\delta^{ij}$ $\in$ $C^{1}([T_1,T_2],H^{\ell}_\delta)$ for
any $\ell \geq 0$ and $0<T_1 < T_2 <T$.
\end{thm}
\begin{proof}
Let $\tilde{\Gamma}^{k}_{ij}$ denote the Christoffel symbols for the
Euclidean Levi-Civita connection on $M\cong \Rbb^n$. Following the
now standard method, see \cite{RF} Sec.~3.3, we first solve the
Hamilton-DeTurck flow
\leqn{LocA4}{
\partial_t g_{ij} = -2R_{ij} + \nabla_i W_i + \nabla_j W_j \;\; , \quad
g(0) = \hat{g} \, , }
where
\leqn{LocA5}{ W_j = g_{jk} W^k := g_{jk} g^{pq}(\Gamma^{k}_{pq} -
\tilde{\Gamma}^{k}_{pq}) \, , }
and $\Gamma^{k}_{ij}$ are the Christoffel symbols for the
Levi-Civita connection derived from $g$. Since $M\cong \Rbb^n$, we
can use global Cartesian coordinates $(x^{1},\ldots,x^n)$ where
$\tilde{\Gamma}^{k}_{ij} = 0$. With respect to the Cartesian
coordinates, the initial value problem \eqref{LocA4} becomes, see
Lemma 2.1 in \cite{Shi},
\lalign{LocA5}{
\partial_t h_{ij} &= g^{ij}\partial_i\partial_j h_{ij}
+ \Half g^{pq}g^{rs}\bigl(\partial_i h_{pr}\partial_j h_{qs}
+ 2\partial_{p}h_{jp}\partial_{q}h_{is}
-  2\partial_{p}h_{jp}\partial_{s}h_{iq} \notag \\
&\qquad \qquad -  2\partial_{j}h_{pr}\partial_{s}h_{iq}
-  2\partial_{i}h_{pr}\partial_{s}h_{jq}
\bigr) \label{LocA5.1}\, , \\
h_{ij}(0) &= \hat{g}_{ij} -\delta_{ij} \in H^k_\delta
\label{LocA5.2} \, , }
where $g_{ij} = \delta_{ij} + h_{ij}$. But $k> n/2+3$ and $\delta
<0$, so we can apply Theorem \ref{locB} to conclude that the
quasi-linear parabolic initial value problem
\eqref{LocA5.1}--\eqref{LocA5.2} has a local solution $h_{ij}(t,x)$
that satisfies
\leqn{LocA6}{ h_{ij}\, , \;g^{ij}-\delta^{ij} \in
C^0([0,T),H^k_\delta) \cap C^{1}([0,T),H^{k-2}_\delta) }
for some $T > 0$,
\leqn{LocA7}{ h_{ij}(t,x)\, ,\; g^{ij}(t,x)\; \in
C^{\infty}((0,T)\times \Rbb^n), }
and $h_{ij} \in C^{1}([T_1,T_2],H^\ell_\delta)$ for any $\ell \geq
0$ and $0<T_1<T_2 < T$. The time-dependent vector field $W^k$ is
given by
\leqn{LocA8}{ W^{k} = g^{ij}\Gamma^{k}_{ij} = \Half
g^{ij}g^{kp}\bigl(\partial_i h_{jp} +\partial_j h_{ip} -\partial_p
h_{ij}\bigr) \, , }
and $W^{k}$ defines a continuous map from $[0,T)$ to
$H^k_\delta(\Rbb^n,\Rbb^n)$ by \eqref{LocA6} and Lemma \ref{SobB}.
Note also that $W^{k}\in C^{\infty}((0,T)\times \Rbb^n)$. Letting
$\psi_t(x)=(\psi^{1}_t(x),\ldots, \psi^{n}_t(x))$ denote the flow of
$W^{k}$ where $\psi_0 = \id$, Theorem \ref{EM1} implies that the
map, shrinking $T$ if necessary, $[0,T) \ni t$ $\mapsto$ $\psi_t$
$\in \mathcal{D}^{k-1}_\delta$ is $C^1$. In particular, this implies
that $\psi_t^i(x) = x^i + \phi^i_t(x)$ where the map $[0,T) \ni t$
$\mapsto$ $\phi_t$ $\in H^{k-1}_\delta$ is $C^1$. But $W^{k}\in
C^{\infty}((0,T)\times \Rbb^n)$, so we also get that $\psi(t,x) \in
C^{\infty}((0,T)\times \Rbb^n)$.

Let $\bar{h}$ denote the pullback of $h$ by the diffeomorphism
$\psi_t$ so that
\leqn{LocA10}{ \bar{h}_{ij} = \bigl(\psi^{*}_t h \bigr)_{ij} =
\bigl(h_{pq}\circ\psi_t\bigr)
\partial_i\psi^p_t\partial_j\psi^q_t \, . }
Then $\bar{h}_{ij} \in C^{0}([0,T),H^{k-2}_\delta)$ by Proposition
\ref{Can3} and Lemma \ref{SobB}. Also, $\bar{h}_{ij}(t,x)\in
C^{\infty}((0,T)\times \Rbb^n)$ by \eqref{LocA7}. Differentiating
\eqref{LocA10} with respect to $t$ yields
\lalign{LocA13}{
\partial_t \bar{h}_{ij} =&
\bigl(\partial_t h_{pq}\circ\psi_t\bigr)\partial_i\psi^p_t\partial_j\psi^q_t
+\bigl(\partial_r h_{pq}\circ\psi_t\bigr)\partial_t\psi_t^r
\partial_i\psi^p_t\partial_j\psi^q_t\notag\\
&+ \bigl(h_{pq}\circ\psi_t\bigr)\bigl(
\partial_i\partial_t\psi^p_t\partial_j\psi^q_t +
\partial_i\psi^p_t\partial_j\partial_t\psi^q_t\bigr) \, . }
Using the same arguments as above, we also find that
$\partial_t\bar{h}_{ij} \in C^{0}([0,T),H^{k-2}_\delta)$.

Finally, let $\bar{g} = \psi_t^*g$. Then $\bar{g}$ is a solution to
the Ricci flow equation, see Ch.~3.3 of \cite{RF}, $\partial_t
\bar{g}_{ij} = -2\bar{R}_{ij}$ with initial data $\bar{g}(0) =
\hat{g}$. Furthermore,
\leqn{LocA14}{ \bar{g}_{ij} -\delta_{ij} = \bar{h}_{ij} +
\delta_{pq}\partial_i\psi^p_t\partial_j\psi^q_t -\delta_{ij} =
\bar{h}_{ij} + \delta_{pq}\partial_i\phi^p_t\partial_j \phi^{q}_t }
and hence $\bar{g}_{ij} -\delta_{ij} \in
C^{1}([0,T),H^{k-2}_\delta)$ since we showed above that
$\partial_j\phi^{i}_t$, $\bar{h}_{ij}$ $\in$
$C^{1}([0,T),H^{k-2}_\delta)$.

Similar arguments show that $\bar{g}_{ij} -\delta_{ij}$,
$\bar{g}^{ij}-\delta^{ij}$ $\in$
$C^{1}([T_1,T_2],H^{\ell}_\delta)$ follows from $h_{ij} \in
C^{1}([T_1,T_2],H^\ell_\delta)$. Also, $\bar{g}_{ij}\in
C^{\infty}((0,T)\times \Rbb^n)$ follows easily from
$h_{ij}(t,x)\, , \psi(t,x) \in C^{\infty}((0,T)\times
\Rbb^n)$.
\end{proof}
\begin{cor}
\label{LocBa} Let $k>n/2+4$ and $g(t)$ be the Ricci flow solution
from Theorem \ref{LocA}. Then $R_{ij} \in
C^1([0,T),H^{k-4}_{\delta-2})$ and $g_{ij}(t) = \hat{g}_{ij} +
f_{ij}(t)$ where $f_{ij} \in C^{1}([0,T),H^{k-4}_{\delta-2})$.
Moreover, if $k>n/2+6$, $\delta < 4-n$ and $\hat{R} \in L^1$ then
$R(t) \in C^{1}([0,T),L^1)$.
\end{cor}
\begin{proof}
Let $h_{ij}= g_{ij}-\delta_{ij}$. Then the  Ricci curvature of $g$
has the form $R_{ij} = B_{ij}(g^{pq},\partial_\ell\partial_m
h_{rs})$ $+C_{ij}(g^{pq},\partial_q h_{rs})$ where $B_{ij}$ and
$C_{ij}$ are analytic functions that are linear and quadratic,
respectively, in their second variables. It follows from the
weighted multiplication Lemma \ref{SobB} that the map
$H^{\ell}_\delta \ni (g^{ij}-\delta^{ij},h_{ij}) \mapsto R_{ij} \in
H^{\ell}_{\delta-2}$ is well defined and analytic for $\eta \leq 0$
and $\ell > n/2$. This proves the first statement.

Integrating $\partial_t g_{ij} = -2 R_{ij}$ with respect to $t$
yields $g_{ij}(t) - \hat{g}_{ij}$ $=$ $-2\int_{0}^{t}R_{ij}(s)ds$.
But $R_{ij} \in C^1([0,T),H^{k-4}_{\delta-2})$, and thus the map
$[0,T) \ni t \mapsto -2\int_{0}^{t}R_{ij}(s)ds \in
H^{k-4}_{\delta-2}$ is well defined and continuously
differentiable. This completes the proof of the second statement.

The Ricci scalar satisfies the equation
\leqn{LocBa5}{
\partial_t R = \Delta R + |\text{Ric}|^2 \, .
}
Integrating this yields $R(t)$ $=$ $\hat{R} +$ $\int_{0}^{t}
\bigl(\Delta R(s) + |\text{Ric}|^2(s)\bigr)ds$ . From Corollary
\ref{LocBa} and the weighted multiplication lemma \ref{SobB}, we see
that $\Delta R + |\text{Ric}|^2 \in
C^{1}([0,T),H^{k-6}_{\delta-4})$. By the weighted H\"{o}lder and
Sobolev inequalities (Lemmata \ref{Holder} and \ref{SobA}), we have
$H^{k-6}_{\delta-4}\subset L^\infty_{\delta-4} \subset L^1$. Thus
$\int_{0}^{t} \bigl(\Delta R(s) + |\text{Ric}|^2(s)\bigr)ds \in L^1$
for all $t\in [0,T)$.
\end{proof}

\begin{rem}
\label{mass} {\rm In  \cite{Bart86} Proposition 4.1, it is
established that the mass of an asymptotically flat metric $g$ of
class $H^k_\delta \subset W^{2,2n/(n-2)}_\delta$  $(k\geq 3)$ is
well defined and given by the formula
\leqn{mass1}{ \text{mass}(g) := \int_{S_\infty}\bigl(\partial_j
g_{ij} -\partial_i g_{jj} \bigr) \, dS^i }
provided $\delta \leq (2-n)/2$ and the Ricci scalar is both
non-negative and integrable. So, by the above corollary and the
maximum principle, see equation (\ref{LocBa5}), an initial
asymptotically flat metric $\hat{g}$ of class $H^k_\delta$, where
$k>n/2+6$ and $\delta < \min\{4-n,(2-n)/2\}$, with non-negative
and integrable Ricci scalar will yield a flow $g(t)$ for which the
Ricci scalar continues to be non-negative and integrable for every
$t>0$. Thus the mass of $g(t)$ remains well defined. Furthermore,
since $g_{ij}-\hat{g}_{ij}\in H^{k-4}_{\delta-2} \subset
W^{1,\infty}_{\delta-2} \subset W^{1,\infty}_{2-n}$, it follows
easily from the definition of the mass that
\leqn{mass2}{ \text{mass}(g(t)) = \text{mass}(\hat{g})\quad
\text{for all $t\geq 0$.} } }
\end{rem}

\begin{thm}
\label{LocB} \mnote{[LocB]}
Suppose $k>n/2+4$, $\delta < 0$, and $\tilde{g}(t)$ and $\bar{g}(t)$
are two solutions to the Ricci flow satisfying $\bar{g}(0) =
\tilde{g}(0)$ and
\eqn{LocB1}{ \tilde{g}_{ij}-\delta_{ij}\, , \;
\tilde{g}^{ij}-\delta^{ij}\, , \; \bar{g}_{ij}-\delta_{ij} \, ,
\bar{g}^{ij} - \delta^{ij} \in C^{1}([0,T),H^k_\delta) . }
Then $\bar{g}(t) = \tilde{g}(t)$ for all $t\in [0,T)$.
\end{thm}
\begin{proof}
Fix $k > n/2+4 $ and $\delta < 0$. To prove uniqueness, we use
Hamilton's method involving harmonic maps \cite{Ham95} as described
in Sec.~3.4 of \cite{RF}. Let $e=\delta_{ij}dx^i dx^j$ denote the
Euclidean metric. As before, $(x^1,\ldots,x^n)$ are Cartesian
coordinates. Given a map $f_{0} : M\cong \Rbb^n \rightarrow M :
x=(x^1,\ldots,x^n) \mapsto (f_0^1(x),\ldots,f_0^n(x))$ and a metric
$g$, the harmonic map flow with respect to the pair $(g,e)$ of
metrics on $M$ is
\leqn{LocB3}{
\partial_t \psi = \Delta_{g,e}\psi \quad : \quad \psi(0) = \psi_0 \;
}
where $\psi_t(x)=(\psi^1_t(x),\ldots,\psi_t^n(x))$ is a time
dependent map from $\Rbb^n$ to $\Rbb^n$ and $\Delta_{g,e}\psi$ is
defined by
\leqn{LocB4}{ (\Delta_{g,e}\psi)^j = g^{pq}\bigl(
\partial_{p}\partial_{q} \psi^j -
\Gamma^r_{pq}\partial_r \psi^j\bigr) \, .
}
As above, $\Gamma^r_{pq}$ are the Christoffel symbols of the
Levi-Civita connection derived from $g$. If we let
\leqn{LocB5}{ \psi^j_t(x) = x^j +\phi^j_t(x) \quad  \text{and} \quad
\psi^j_0(x) = x^j +  \phi^j_0(x)\, , }
then we can write \eqref{LocB3} as
\leqn{LocB6}{
\partial_t \phi^j = g^{pq}\bigl(
\partial_{p}\partial_{q} \phi^j -
\Gamma^r_{pq}\partial_r \phi^j -\Gamma^j_{pq}\bigr) \quad , \quad
\phi^j(0) = \phi^j_0 \, . }
Suppose $g$ is a time dependent metric that satisfies $g \in
C^{0}([0,T),H^{k}_\delta)$. Then the continuity of the
differentiation operator and Lemma \ref{SobB} imply that
$\Gamma^r_{pq} \in C^1([0,1),H^{k-1}_\delta)$. So if $\phi_0^j \in
H^{k-1}_{\delta}$, then there exists a unique solution $\phi^{j}$
$\in C^{0}([0,T),H^{k-1}_\delta)$$\cap$ $C^1([0,T),H^{k-3}_\delta)$
to \eqref{LocB6} by Theorem \ref{locB}. If $\psi_0 \in
\mathcal{D}^{k-1}_\delta$, then Theorem \ref{Can3A} implies,
shrinking $T$ if necessary, that $\psi$ $\in
C^{0}([0,T),\mathcal{D}^{k-1}_\delta)$ $\cap$ $
C^1([0,T),\mathcal{D}^{k-3}_\delta)$.

Suppose $\tilde{g}(t)$ and $\bar{g}(t)$ are two solutions to the
Ricci flow such that
\eqn{LocB9}{ \tilde{g}_{ij}-\delta_{ij}\, , \;
\tilde{g}^{ij}-\delta^{ij}\, , \; \bar{g}_{ij}-\delta_{ij} \, ,
\bar{g}^{ij} - \delta^{ij} \in C^{1}([0,T),H^k_\delta) \, }
and $\bar{g}(0) = \tilde{g}(0)$. Let $\tilde{\psi}$, $\bar{\psi}$
$\in C^{0}([0,T),\mathcal{D}^{k-1}_\delta)$ $\cap
C^1([0,T),\mathcal{D}^{k-3}_\delta)$ be solutions to the harmonic
map flow with respect to the metric pairs $(\tilde{g},e)$ and
$(\bar{g},e)$  with initial conditions
$\tilde{\psi}(0)=\bar{\psi}(0) = \id$. Letting $\tilde{h}_{ij}$ $:=
(\tilde{\psi}_*\tilde{g})_{ij}-\delta_{ij}$ and $\bar{h}_{ij}$ $:=
(\bar{\psi}_*\bar{g})_{ij}-\delta_{ij}$, the same arguments as in
the proof of Theorem \ref{LocA} show that
\eqn{LocB13}{ \tilde{h}_{ij}\, , \;
(\tilde{\psi}_*\tilde{g})^{ij}-\delta^{ij} \, , \; \bar{h}_{ij} \, ,
(\bar{\psi}_*\bar{g})^{ij}-\delta^{ij} \; \in
C^{0}([0,T),H^{k-2}_\delta)\cap C^{1}([0,T),H^{k-4}_\delta) \, . }
But $\tilde{\psi}_*\tilde{g}$ and $\bar{\psi}_*\bar{g}$ both satisfy
the Hamilton-DeTurck flow \eqref{LocA4} (see Sec.~3.4.4 of
\cite{RF}) or equivalently $\bar{h}_{ij}$ and $\tilde{h}_{ij}$ both
satisfy the parabolic equation \eqref{LocA5.1} with initial
condition $\bar{h}_{ij}(0) = \tilde{h}_{ij}(0)$. By uniqueness of
solutions to this equation (see Theorem \ref{locB}) we must have
$\bar{h}_{ij}(t) = \tilde{h}_{ij}(t)$ or equivalently
$(\tilde{\psi}_*\tilde{g})(t)= (\bar{\psi}_*\bar{g})(t)$ for all
$t\in [0,T)$. So by Lemma 3.27 in \cite{RF}, the time dependent
diffeomorphisms $\tilde{\psi}_t$ and $\tilde{\psi}_t$ are flows for
the time dependent differential equation $dx^j/dt$ $=$ $W^j(t,x)$
that satisfy $\tilde{\psi}_0$ $=$ $\bar{\psi}_0$ $=$ $\id$. Here,
$W^j$ is the vector field defined by $W^j$ $= g^{pq}\Gamma^j_{pg}$.
By standard uniqueness theorems for solutions to ordinary
differential equations, we can conclude that $\tilde{\psi}_t$ $=$
$\bar{\psi}_t$ for all $t$ $\in$ $[0,T)$. It follows that
$\tilde{g}(t)$ $=$ $\bar{g}(t)$ for all $t$ $\in$ $[0,T)$ and the
proof is complete.
\end{proof}

\subsection{A Continuation Principle}

\noindent The following theorem shows that if local existence in
time fails to extend indefinitely to give global future existence,
then curvature diverges in finite time.
\begin{thm}
\label{cont} \mnote{[cont]} Suppose $k>n/2 +4$, $\delta < 0$ and
$\hat{g}$ is an asymptotically flat metric of class $H^k_\delta$.
Then Ricci flow $\partial_t g_{ij} = -2R_{ij}$ with the initial
condition $g(0)=\hat{g}$ has a unique solution on a maximal time
interval $0\leq t < T_M \leq \infty$. If $T_M < \infty$ then
\leqn{cont1} {\limsup_{t\to T_M}
%\underset{\;t\nearrow T_M} {\overline{\lim}\;\;}
\sup_{x\in\Rbb^n}|{\rm Rm}(t,x)|_{g(t,x)} = \infty \, . } Moreover,
for any $T \in [0,T_M)$, $K= \sup_{0\leq t\leq T}
\sup_{x\in\Rbb^n}|{\rm Rm}(t,x)|_{g(t,x)} < \infty$ and
\leqn{cont2}{ e^{-2KT}\hat{g} \leq g(t) \leq e^{2KT}\hat{g} \quad
\text{for all $t\in [0,T]$.} }
\end{thm}
\begin{proof}
For $\hat{g}\in H^k_\delta$ with $k>n/2+4$ and $\delta < 0$, let
$[0,T_M)$ be the maximal time interval of existence for a solution
$g(t)$ to Ricci flow. Suppose that $T_M < \infty$  and that
\leqn{cont3}{ K := \sup_{0\leq t< T_M} \sup_{x\in\Rbb^n}|{\rm
Rm}(t,x)|_{g(t,x)} < \infty. }
For each $t\in [0,T_M)$, the metric is $g(t)$ is asymptotically flat
and hence $g(t)$ is a solution to Ricci flow for which the maximum
principle holds. It follows that Proposition 6.48 of \cite{RF}
applies.
%is valid.
So for each $m\in \Nbb_0$, there exists a constant $c_m$ such that
\leqn{cont4}{ |D^I g_{ij}(t,x)|+ |D^I g^{ij}| \leq c_m \quad
\text{for all $|I|= m$ and $(t,x)\in [0,T_M)\times\Rbb^n$}, }
where $g_{ij}$ are the metric components in Cartesian coordinates
and $D^I =
\partial^{I_1}_1\ldots \partial^{I_n}_n$.

{}From the proof of Theorem \ref{LocB}, we get that for each
$\ttl\in [0,T_M)$ there exists an interval $I_{\ttl} :=
[\ttl,T_{\ttl}) \subset [0,T_M)$ and a map $\psi_{\ttl}(t,x) =
(\psi_{\ttl}^1(t,x),\ldots,\psi_{\ttl}^n(t,x))$ of $\Rbb^n$ to
$\Rbb^n$ such that
\eqn{cont5}{ \psi_{\ttl}\in C^{0}(I_{\ttl},\Dc^{k-1}_{\delta})\cap
C^{1}(I_{\ttl},\Dc^{k-3}_{\delta})}
and $\psi_{\tilde{t}}$ satisfies harmonic map flow (i.e.
\eqref{LocB3}) with initial condition $\psi_{\ttl}=\id$. Since
$\psi_{\tilde{t}}$ satisfies a linear equation, see \eqref{LocB3},
$\psi_{\tilde{t}}$ will continue to exist as long as $g(t)$ does.
Thus we can solve \eqref{LocB3} on the interval $[\ttl,T_M)$
although it may fail to define a diffeomorphism for some time less
than $T_M$. Also, the metric $\gt(t):= (\psi_{\ttl})_*g(t)$
satisfies
\eqn{cont6}{ \gt^{ij} -\delta^{ij}\, , \gt_{ij}-\delta_{ij} \in
C^{0}(I_{\ttl},H^{k-2}_\delta)\cap C^{1}(I_{\ttl},H^{k-4}_\delta) }
and $\htl_{ij} := \gt_{ij}-\delta_{ij}$ is a solution of
Hamilton-DeTurck flow \eqref{LocA5.1} on the time interval
$I_{\ttl}$ with initial condition
$\htl_{ij}(\ttl)=g_{ij}(\ttl)-\delta_{ij}$.

We now use the harmonic map flow equation \eqref{LocB3} to derive
$C^k_b$ bounds on $\psi_{\ttl}$ to estimate the length of time for
which $\psi_{\ttl}$ remains a diffeomorphism. Let $\phi^j(t,x) =
\psi^j_{\ttl}(t,x)-x^j$ and define $|\phi|^2 =
\delta_{ij}\phi^i\phi^j$. Then from \eqref{LocB3}, or equivalently
\eqref{LocB6}, $|\phi|^2$ satisfies
\leqn{cont7}{
\partial_t|\phi|^2 = g^{pq}\partial_p\partial_q |\phi|^2
-\delta_{ij}g^{pq}\partial_p\phi^i\partial_q\phi^j
-g^{pq}\Gamma^r_{pq}\partial_r|\phi|^2 -\Half\delta_{ij}\phi^i
g^{pq}\Gamma^j_{pq} \, . }
So by \eqref{cont4}, there exists a constant $C$ independent  of
$\ttl$  such that
\leqn{cont8}{
\partial_t|\phi|^2 - g^{pq}\partial_p\partial_q |\phi|^2
+g^{pq}\Gamma^r_{pq}\partial_r|\phi|^2 \leq |\phi|^2 + C\, . }
Since $\lim_{|x|\rightarrow \infty}|\phi|^2(t,x) = 0$ for all $t\in
[\ttl,T_M)$ and $|\phi|^2(0,x) = 0$, we get via the maximum
principle, see Theorem 4.4 in \cite{RF}, that
\leqn{cont9}{ |\phi|^2(t,x) \leq C(\exp(t-\ttl)-1) \quad \text{ for
all $(t,x) \in [\ttl,T_M)\times \Rbb^n$.} }
Next, differentiating \eqref{LocB3} we find that
\lalign{cont10}{
\partial_t|D\psi_{\ttl}|^2 = g^{pq}&\partial_p\partial_q |D\psi_{\ttl}|^2
-2g^{pq}\delta_{jl}\delta^{ki}\partial_q
\partial_k\psi_{\ttl}^l\partial_p \partial_i\psi_{\ttl}^j
-g^{pq}\Gamma^r_{pq}\partial_r|D\psi_{\ttl}|^2 \notag \\
&+2\partial_i g^{pq}\delta_{jl}\delta^{ki}
\partial_k\psi_{\ttl}^l\partial_p \partial_q\psi_{\ttl}^j - 2
\partial_i(g^{pq}\Gamma^r_{pq})\delta_{jl}\delta^{ki}
\partial_k\psi_{\ttl}^l \partial_r\psi_{\ttl}^j.
}
where $|D\psi_{\ttl}|^2 :=
\delta_{ij}\delta^{kl}\partial_k\psi^i_{\ttl}\partial_l\psi^j_{\ttl}$.
Using \eqref{cont4}, we obtain the inequalities
\lalign{cont11}{ &
\partial_t|D\psi_{\ttl}|^2 -
g^{pq}\partial_p\partial_q|D\psi_{\ttl}|^2 +
g^{pq}\Gamma^r_{pq}\partial_r|D\psi_{\ttl}|^2\notag\\
& \leq -2g^{pq}\delta_{jl}\delta^{ki}\partial_q
\partial_k\psi_{\ttl}^l
\partial_p \partial_i D\psi^j_{\ttl}
+ \epsilon|D\psi_{\ttl}|^2 + C_1(1+1/\epsilon)|D\psi_{\ttl}|^2
\quad (\epsilon
> 0) }
and
\leqn{cont12}{ -2g^{pq}\delta_{jl}\delta^{ki}\partial_q
\partial_k\psi_{\ttl}^l
\partial_p \partial_i\psi_{\ttl}^j \leq - C_2 |D\psi_{\ttl}|^2 \, . }
for some constants $C_1$ and $C_2$ that are independent of $\epsilon
> 0 $, $\ttl$ and $t\in [\ttl,T_M)$. Setting $\ep = C_2$ yields
\lalign{cont13}{
\partial_t|D\psi_{\ttl}|^2 - g^{pq}\partial_p\partial_q|D\psi_{\ttl}|^2 +
g^{pq}\Gamma^r_{pq}\partial_r|D\psi_{\ttl}|^2 \leq
C_1(1+1/C_2)|D\psi_{\ttl}|^2 \, . }
Since $\lim_{|x|\rightarrow \infty }|D\psi_{\ttl}|^2(t,x) = n$ for
all $t\in [\ttl,T_M)$ and $|D\psi_{\ttl}|^2(\ttl,x)=n$, the maximum
principle implies that there exists a constant $C$ independent of
$\ttl$ for which the following estimate holds
\leqn{cont14}{ ||D\psi_{\ttl}|^2(t,x)-n| \leq C\exp((t-\ttl)-1)
\quad \text{for all $(t,x)\in [\ttl,T_M)\times \Rbb^n$.} }
Differentiating \eqref{LocB3} again and letting $|D^2\psi_{\ttl}|
=
\delta_{ij}\delta^{kl}\delta^{pq}\partial_{kp}\psi^i_{\ttl}
\partial_{lq}\psi^j_{\ttl}$,
we find, using similar arguments, that there exists a constant
$C>0$ independent of $\ttl$ such that
\leqn{cont15}{|D^2\psi_{\ttl}|^2(t,x) \leq \exp(C(t-\ttl)) \quad
\text{for all $(t,x)\in [\ttl,T_M)\times \Rbb^n$.}}

Let $J(\psi_{\ttl})=\det(\partial_j\phi^i_{\ttl})$ denote the
Jacobian of the map $\psi_{\ttl}$. Since $J(\psi_{\ttl})=1$, the
estimates \eqref{cont9} and \eqref{cont14} show that there exists a
$\bar{t}\in (0,T_M)$ and a constant $C>1$ such that
\leqn{cont16}{ 0<1/C \leq J(\psi_{\ttl})(t,x) \leq C  \quad
\text{for all $(t,x) \in [\bar{t},T_M)\times \Rbb^n$.} }
Combining this estimate with \eqref{cont9}, \eqref{cont14}, and
\eqref{cont15}, we have
\leqn{cont17}{ |D^I(\psi^{-1}_{\ttl}(t,x)-x)|\leq C \quad \text{for
all $|I|\leq 2$ and $(t,x)\in [\bar{t},T_M)\times \Rbb^n$.}}
Notice that this estimate along with Lemma \ref{Can3B} shows that
$I_{\bar{t}}=[\bar{t},T_M)$ and that
\leqn{cont18}{ |D^I \tilde{h}_{ij}(t,x)| + |g^{ij}(t,x)| \leq C
\quad \text{for $|I|\leq 1$ and all $(t,x)\in [\bar{t},T_M)$.}}
But $\tilde{h}$ satisfies \eqref{LocA5.1}, and so the estimate
\eqref{cont18} and the continuation principle of Theorem \ref{locB}
imply that there exists a $T>T_M$ such that $\tilde{h}_{ij}(t,x)$
extends to a solution on $[\bar{t},T)\times \Rbb^n$ of the class
\leqn{cont19}{ \tilde{h}_{ij}=\tilde{g}_{ij}-\delta_{ij} \,
,\gt^{ij} -\delta^{ij}\, \in C^{0}([\bar{t},T),H^{k-2}_\delta)\cap
C^{1}([\bar{T},T), H^{k-4}_\delta) \, . }
By the proof of Theorem \ref{LocA} and \ref{LocB}, $\tilde{h}_{ij}$
produces a unique solution to Ricci flow satisfying
$g_{ij}-\delta_{ij}$, $g^{ij}-\delta^{ij}$$\in$
$C^{1}([\bar{t},T),H^{k-4}_\delta) $ and $\bar{g}(\bar{t})$ $=$
$g(\bar{t})$. Thus $\bar{g}(t) = g(t)$ for all $t\in [0,T_M)$. Since
$T>T_M$ this contradicts $T_M$ being the maximal existence time. So
we must either have $T_M=\infty$ or $\limsup_{t\nearrow T_M}$
$\sup_{x\in\Rbb^n}|{\rm Rm}(t,x)|_{g(t,x)}$ $<$ $\infty$.  This
proves the first statement. The second statement follows from a
straightforward adaptation of Corollary 6.50 in \cite{RF}.
\end{proof}

We note that as in the compact case the continuation criterion
(\ref{cont1}) can be strengthened to $\lim_{t\nearrow T_M}$
$\sup_{x\in\Rbb^n}|{\rm Rm}(t,x)|_{g(t,x)}$ $=$ $\infty$ but we will
not pursue this here.

\section{Rotational Symmetry}
\setcounter{equation}{0}

\subsection{The Coordinate System}

\noindent We now restrict our attention to flows evolving from a
fixed initial initial metric that (i) is rotationally symmetric and
admits no minimal hyperspheres, and
(ii) is asymptotically flat of class $H^k_\delta$ with $\delta <0$
and $k>\frac{n}{2}+4$.
%(ii) lies in $H^k_\delta$ with
%$\delta < \min\{4-n,1-n/2\}$ and $k > n/2 +6$, and (iii) has
%$\hat{R} \in L^1$ whenever $\hat{R}\geq 0$.
%In this section we use $T_M$ to denote the maximal existence time
%for the flow, from Theorem \ref{cont}.
%
%We can always choose global coordinates in which such a metric can
%be expressed as
%
%
%\leqn{eq4.1}{ \hat{g} = a^2(r)dr^2 + r^2 g_{\rm can} }
%
%
%where $r = |x|:=\sqrt{\delta_{ij}x^i x^j}$ and $g_{\rm can}$ is the
%canonical round metric on ${\mathbb S}^{n-1}$.
In an attempt to manage the several constants that will appear from
here onward, we will sometimes use the notation $C^+_x$ to denote a
constant that bounds a quantity $x$ from above; dually, $C^-_x$ will
sometimes be used to denote a constant that bounds $x$ from below.

\begin{rem}\label{locrem} \mnote{[locrem]}
$\;$ {\rm
\begin{enumerate}
\item[(i)] By Theorem \ref{LocA}, there exists a solution $\bar{g}(t)$ to
Ricci flow satisfying
\leqn{eq4.1}{ \bar{g}_{ij}-\delta_{ij},\; \bar{g}^{ij}-\delta^{ij}
\in C^1([0,T_M),H^{k-2}_{\delta}) \, , \quad \bar{g}(t,x) \in
C^\infty((0,T_M)\times \Rbb^n),  }
and $\bar{g}(0)=\hat{g}$. \item[(ii)] From (\ref{eq4.1}) and the
weighted Sobolev embedding (see Lemma \ref{SobA}), it follows that
$\bar{g}(t) \in C^1([0,T_M),C^2_{\delta})$ and hence there exists
a time dependent constant $C(t)$ such that
\leqn{eq4.2}{ |D_x^I \bar{g}_{ij}(t,x)| \leq
\frac{C(t)}{(1+|x|^2)^{(|\delta|+|I|)/2}} }
for all $(t,x)$ $\in$ $[0,T_M)\times \Rbb^n$, and $|I|\leq 2$.
\item[(iii)] Since Ricci flow preserves isometries, each metric
$g(t)$ is rotationally symmetric and hence
\leqn{eq4.3}{ \bar{g}(t,x) = q^2(t,r) dr^2 + h^2(t,r) g_{\rm can} }
for functions $q(t,r)$ and $h(t,r)$ which are $C^1$ in $t$, $C^2$
in $r$, $C^\infty$ in $t$ and $r$ for $t
> 0$, and satisfy
\lgath{}{ q(0,r) = a(r) \, ,\quad  \quad
h(0,r) = r\,  , \label{eq4.4} \\
|\partial_r^s(q^2(t,r)-1)| \leq \frac{C(t)}{(1+r)^{|\delta|+s}}
\quad s=0,1,2 \, , \label{eq4.5} \\
|\partial_r^s(r^{-2} h^2(t,r)-1)| \leq
\frac{C(t)}{(1+r)^{|\delta|+s}} \quad s=0,1,2 \, . \label{eq4.6} }
%
%
%\item[(iv)] As discussed in remark \ref{mass}, if $\hat{R}\geq 0$,
%then the restrictions on $\delta$ and $k$ guarantee that the mass of
%both the initial metric $\hat{g}$ and the evolved metric
%$\bar{g}(t)$ are well-defined and satisfy ${\rm
%mass}(\bar{g}(t))={\rm mass}(\hat{g})$ for all $t\in [0,T_M)$.
\end{enumerate}
}
\end{rem}

Since  $\partial_r h(0,r) = \partial_r r  =1$, it follows that there
exist constants $0<C_{\partial_r h}^- \leq 1$, $C_{\partial_r
h}^+\geq 1$, such that
\leqn{eq4.7}{ 0 < C_{\partial_r h}^- \leq  \partial_r h(t,r) \leq
C_{\partial_r h}^+ \quad \text{for all $(t,r) \in [0,T]\times
(0,\infty)$} }
for some $T>0$. Note that $T$ has no {\it a priori} relation to
$T_M$, the maximal existence time of the flowing metric
(\ref{eq4.3}). However, let ${\tilde T}$ be the largest time such
that (\ref{eq4.7}) holds whenever $T<{\tilde T}$. We will show in
Subsection 4.3 that we can take ${\tilde T}=T_M$.

Letting $(\theta^A)$ denote  angular coordinates on the sphere
${\mathbb S}^{n-1}$, the map
\leqn{eq4.8}{ \psi_t(r,\theta^A) = (h(t,r),\theta^A) }
defines a $C^2$ diffeomorphism on $\Rbb^n$  for each $t \in
[0,{\tilde T})$ which is smooth for all $t>0$. So then
\leqn{eq4.9}{ \psi^{-1}_t(r,\theta^A) = (\rho(t,r),\theta^A) }
for a function $\rho(t,r)$ that is $C^1$ in $t$, $C^2$ in $r$,
$C^\infty$ in $r$ and $t$ for $t>0$, and satisfies
\leqn{eq4.10}{ h(t,\rho(t,r))=r \, , \quad \rho(t,h(t,r)) = r\, ,
\quad \text{and} \quad \rho(0,r) =r }
for all  $(t,r) \in [0,T]\times (0,\infty)$. Next, define
\leqn{eq4.11}{ g(t) := (\psi_t^{-1})^* \bar{g}(t). }
Then we finally obtain that
\leqn{eq4.12}{ g(t) = f^2(t,r) dr^2 + r^2 g_{\rm can} }
where
\leqn{eq4.13}{ f(t,r) = \frac{q(t,\rho(t,r))}{\partial_r
h(t,\rho(t,r))} \quad \text{for all $(t,r) \in [0,{\tilde T})\times
(0,\infty)$.} }

Note that $f(t,r)$ is $C^1$ in $t$, is $C^2$ in $r$, and $C^\infty$
in $r$ and $t$ for $t>0$. As well,
\begin{equation}
\lim_{r\to\infty}f^2(t,r)=1 \label{eq4.14}
\end{equation}
(proof: from (\ref{eq4.5}) we have $q^2\to 1$ and from (\ref{eq4.6})
it's easy to check that $\partial_r h\to 1$; then apply these in
(\ref{eq4.13})). Finally note that the mean curvature of
constant-$r$ hyperspheres is
\begin{equation}
H=\frac{1}{rf}\ , \label{eq4.15}
\end{equation}
so a minimal hypersphere occurs iff $f$ diverges at finite $r$ and
some $t\in[0,{\tilde T}]$. We show in the following subsection that
such a divergence cannot develop.

\subsection{Ricci Flow in Area Radius Coordinates}

\noindent The metric (\ref{eq4.12}) is a solution of the
Hamilton-DeTurck flow (\ref{eq1.2}), at least for $t\in[0,{\tilde
T})$. Now from (\ref{eq4.12}) we can directly compute the Ricci
curvature and obtain
\begin{equation}
{\rm Ric} = \frac{(n-1)}{rf(t,r)} \frac{\partial f}{\partial r} dr^2
+ \left [ (n-2) \left ( 1 - \frac{1}{f^2(t,r)} \right )
+\frac{r}{f^3(t,r)}\frac{\partial f}{\partial r}\right ] g_{\rm can}
\ . \label{eq4.16}
\end{equation}
We can then use the components of the flow equation (\ref{eq1.2})
normal to $\frac{\partial}{\partial r}$ to determine $\xi$,
expressed as a 1-form, to be $\xi=\xi_1(t,r)dr$ where
\begin{equation}
\xi_1=\left [ \frac{(n-2)}{r} \left ( f^2(t,r) - 1 \right )
+\frac{\frac{\partial f}{\partial r}}{f(t,r)}\right ] \ .
\label{eq4.17}
\end{equation}
We can then write the $rr$-component of (\ref{eq1.2}) as a
differential equation for $f$ and use (\ref{eq4.17}) to eliminate
$\xi$ from this equation. The result is
\begin{eqnarray}
\frac{\partial f}{\partial t}&=&\frac{1}{f^2}\frac{\partial^2
f}{\partial r^2} -\frac{2}{f^3}\left ( \frac{\partial f}{\partial r}
\right )^2 + \left ( \frac{(n-2)}{r}-\frac{1}{rf^2} \right )
\frac{\partial f}{\partial r}\nonumber\\
&&-\frac{(n-2)}{r^2f}\left ( f^2 -1 \right )\ . \label{eq4.18}
\end{eqnarray}
This is our master equation upon which our global existence proof is
based. Obviously $f(t,r)=1$ (flat space) is a solution, as is
$f=const\neq 1$ when $n=2$ (flat cone) but not for $n>2$ .

We will now prove that minimal hyperspheres cannot form along the
flow if none are present initially. A variant of this argument will
be employed several times over in Section 5. Our technique is to
prescribe limits as $r\to\infty$ and as $r\to 0$ on $f(t,r)$ or,
depending on the situation, an expression involving $f$ (and, in the
next section, its radial derivative as well). These limits
constitute time-dependent bounds on the behaviour of the geometry
over the time interval $[0,{\tilde T})$. But if the flow exists
subject to these limits, then maximum principles will give bounds
expressed solely in terms of the initial conditions. The bounds are
therefore uniform in time and independent of ${\tilde T}$.

To see how this works, express (\ref{eq4.18}) in terms of the
variable
\begin{equation}
w(t,r):=f^2(t,r)-1\ . \label{eq4.19}
\end{equation}
Then, working from (\ref{eq4.18}), we see that $w$ obeys
\begin{equation}
\frac{\partial w}{\partial t}= \frac{1}{f^2} \frac{\partial^2
w}{\partial r^2} - \frac{3}{2f^4} \left [ \frac{\partial w}{\partial
r} \right ]^2 +\left [ \frac{n-2}{r}-\frac{1}{rf^2} \right ]
\frac{\partial w}{\partial r}-\frac{2(n-2)}{r^2}w\ . \label{eq4.20}
\end{equation}

%The solution $f(t,r)>0$ of (\ref{eq4.18})
%is suitably behaved at $r=0$ and infinity for $t\in [0,{\tilde T})$;
%specifically,
Since $f(t,r)$ solves (ref{eq4.18}) and obeys $\lim_{r\to 0}
f^2(t,r)=1=\lim_{r\to \infty} f^2(t,r)$, the corresponding $w=f^2-1$
will solve (\ref{eq4.20}) with $\lim_{r\to 0} w(t,r)=0=\lim_{r\to
\infty} w(t,r)$.

\begin{prop}\label{Prop4.2}
Suppose that $w(t,r)$ is a classical solution of (\ref{eq4.20}) for
$(t,r)\in [0,{\tilde T})\times[0,\infty)=:{\tilde D}$ and that
$\lim_{r\to 0} w(t,r)=0 =\lim_{r\to \infty} w(t,r)$ for all
$t\in[0,{\tilde T})$. Then there exist constants $C^-_w\le 0$ and
$C^+_w\ge 0$ such that $C^-_w\le w(t,r)\le C^+_w$ for all $(t,r)\in
{\tilde D}$.
\end{prop}

\begin{proof}
First choose positive constants $0<r_1<r_2$ and restrict the domain
to $r\in[r_1,r_2]$. Let $T<{\tilde T}$. By the maximum principle, if
the maximum of $w$ on $[0,T]\times [r_1,r_2]$ is positive, it must
lie on the {\it parabolic boundary} $P$ (which consists of those
points where either $t=0$, $r=r_1$, or $r=r_2$). But now take the
limits $r_1\to 0$ and $r_2\to\infty$. By assumption, $w(t,r_1)$ and
$w(t,r_2)$ tend to zero in these limits, so for $r_1$ small enough
and $r_2$ large enough, the maximum, if it is positive, lies on the
{\it initial boundary} $\{ (t,r)\vert t=0\}$ (and since $w(0,0)=0$,
even when the maximum is zero it is realized on the initial
boundary). Finally, take $T\to {\tilde T}$. This proves
\begin{equation}
C^+_w:=\max_{r\in[0,\infty)} \{w(0,r)\}=\max_{\tilde D} \{ w(t,r) \}
\ge 0\ . \label{eq4.21}
\end{equation}

Dually, by the minimum principle, if the minimum of $w$ on
$[0,T]\times [r_1,r_2]$ is negative, it must lie on $P$, and the
argument proceeds as before, yielding
\begin{equation}
C^-_w:=\min_{r\in[0,\infty)} \{w(0,r)\}=\min_{\tilde D} \{ w(t,r)
\}\le 0\ . \label{eq4.22}
\end{equation}
%\qed
\end{proof}

\begin{cor}\label{Cor4.3}
Define constants $C^{\pm}_{f^2}$ such that
$0<C^-_{f^2}:=\min_{r\in[0,\infty)} \{a^2(r)\}$ and let
$C^+_{f^2}:=\max_{r\in[0,\infty)} \{a^2(r)\}$ ($a(r)$ is defined in
(\ref{eq4.4}). Then
\begin{equation}
0<C^-_{f^2}\le f^2(t,r)\le C^+_{f^2}\ . \label{eq4.23}
\end{equation}
for all $(t,r)\in {\tilde D}=[0,{\tilde T})\times [0,\infty)$.
\end{cor}

\begin{proof}
Using $w:=f^2-1$ and noting in particular that
$w(0,r)=f^2(0,r)-1=a^2(r)-1$, apply Proposition (\ref{Prop4.2}) and
use $C^{\pm}_w +1=C^{\pm}_{f^2}$.
\end{proof}

Now we say that a minimal hypersphere forms along the flow iff
$f(t,r)$ diverges in ${\tilde D}=[0,{\tilde T})\times [0,\infty)$.

\begin{cor}\label{Cor4.4}
If no minimal sphere is present initially then none
forms.
\end{cor}

\begin{proof}
{}From Corollary \ref{Cor4.3}, the classical solutions $f$ of
(\ref{eq4.18}) developing from initial data (\ref{eq4.1}) are
bounded uniformly in $t$ on $[0,{\tilde T})$.
\end{proof}

\subsection{The Continuation Principle in Area Radius Coordinates}

\noindent To adapt the continuation principle of Section 3.2 to the
rotationally symmetric case, we must deal with the following point.
While we can assume the solution of Ricci flow in the coordinate
system (\ref{eq4.3}) to exist for all $t<T_M$, the diffeomorphism
transforming the coordinates to those of (\ref{eq4.12}) is, so far,
only defined for $t<{\tilde T}$, and perhaps ${\tilde T}<T_M$.

\begin{prop}\label{Prop4.5}
${\tilde T}= T_M$.
\end{prop}

\begin{proof}
Let $K= \sup_{0\leq t \leq T'}\norm{{\rm Rm}}_{L^\infty}$. But
$\bar{R}_{ijkl}$ is bounded on $[0,T']$ (indeed, on any closed
subinterval of $[0,T_M)$), so we can use (\ref{cont2}), which states
that for all $(t,r)\in [0,T']\times [0,\infty)$
\begin{equation}
e^{-2KT'}C^-_{f^2}\leq e^{-2KT'}a^2(r)\leq q^2(t,r) \leq
e^{2KT'}a^2(r) \leq e^{2KT'}C^+_{f^2}\ . \label{eq4.24}
\end{equation}
Here the inner two inequalities come from (\ref{cont2}) and the
outer two are just the definitions of the constants $C^{\pm}_{f^2}$.

Assume by way of contradiction that ${\tilde T}<T_M$. If we restrict
attention to $t\in[0,{\tilde T})$ then we can divide (\ref{eq4.24})
by (\ref{eq4.23}). This yields
\begin{equation}
0<e^{-2KT'}\frac{C^-_{f^2}}{C^+_{f^2}}\le
\frac{q^2(t,r)}{f^2(t,r)}\le e^{2KT'}\frac{C^+_{f^2}}{C^-_{f^2}}
\label{eq4.25}
\end{equation}
on $[0,{\tilde T})$. Using (\ref{eq4.13}), we can rewrite this as
\begin{equation}
0<e^{-2KT'}\frac{C^-_{f^2}}{C^+_{f^2}}\le \frac{\partial h}{\partial
r} \le e^{2KT'}\frac{C^+_{f^2}}{C^-_{f^2}} \label{eq4.26}
\end{equation}
on $[0,{\tilde T})$. We see by comparison of this to (\ref{eq4.7})
that the constants that appear in (\ref{eq4.7}) are independent of
$T$. But the $\le$ signs give closed relations so, by relaxing the
constant bounds slightly if necessary (keeping the lower bound
positive of course), we can extend (\ref{eq4.26}) (equivalently,
(\ref{eq4.7})) to $[0,{\tilde T}]$ and then to some interval
$[0,T')\supset [0,{\tilde T}]$. This contradicts the assumption that
${\tilde T}<T_M$, and since necessarily ${\tilde T} \le T_M$ we must
therefore conclude that ${\tilde T}= T_M$.
\end{proof}
\noindent Thus the diffeomorphism (\ref{eq4.8}--\ref{eq4.11}) is
defined for all $t\in[0,T_M)$. The square of the norm of the
curvature tensor is given by
\leqn{eq4.27}{ |{\rm Rm}|^2 = R_{ijkl}R^{ijkl} = 2(n-1)\lambda^2_1 +
(n-1)(n-2)\lambda_2^2 }
where
\leqn{eq4.28}{ \lambda_1(t,r) = \frac{1}{rf^3(t,r)}\frac{\partial
f(t,r)}{\partial r} }
and
\leqn{eq4.29}{ \lambda_2(t,r) = \frac{1}{r^2} \left (
1-\frac{1}{f^2(t,r)}\right ) }
are the sectional curvatures in planes containing and orthogonal to
$dr$, respectively. Now in terms of the curvature tensor
$\bar{R}_{ijkl}$ of ${\bar g}(t)$ we have that
\leqn{eq4.30}{ |{\rm Rm}| = |\overline{{\rm Rm}}|\circ\psi^{-1}_t \,
. }
But $\bar{R}_{ijkl}$ is bounded on any interval $[0,T']$ with
$T'<T_M$ and thus the sectional curvatures are bounded functions of
$(t,r)\in [0,T']\times[0,\infty)$, using Proposition 4.5. Thus
\begin{eqnarray}
C^-_{\lambda_1}(t)\le\lambda_1(t,r)&=& \frac{1}{rf^3}\frac{\partial
f}{\partial r} \le C^+_{\lambda_1}(t)\ ,
\label{eq4.31} \\
C^-_{\lambda_2}(t)\le\lambda_2(t,r) &=& \frac{1}{r^2}\left (
1-\frac{1}{f^2(t,r)}\right ) \le C^+_{\lambda_2}(t) \ ,
\label{eq4.32}
\end{eqnarray}
for all $t\in [0,T_M)$. In particular, the limits $r\to 0$ of these
quantities exist at each fixed $t$. It also follows easily from the
fall-offs (\ref{eq4.5}, \ref{eq4.6}) that
\leqn{eq4.33}{ \lim_{r\rightarrow \infty}
r^{-|\delta|-s}\partial_r(f^2(t,r)-1) = 0}
for $s=0,1,2$ and all $t\in [0,T')$. Thus we have that

%\begin{lem}
%\label{fpde} \mnote{[fpde]}

\begin{prop}\label{Prop4.6}
The function $f(t,r)$ given by (\ref{eq4.13}) solves the PDE
(\ref{eq4.18}) on the region $[0,T_M)\times (0,\infty)$, equals
$a(r)$ at time $t=0$, and satisfies the boundary conditions
\leqn{eq4.34}{  \lim_{r\to 0} \frac{1-f^2(t,r)}{r^2} = L_1(t) \,  ,
\quad \lim_{r\to 0} \frac{\partial_r f(t,r)}{r} = L_2(t)  \, ,}
for locally bounded functions $L_1 , L_2 : [0,T_M)\to \mathbb R$ and
%fpde2
\leqn{eq4.35}{ \lim_{r\to \infty}
r^{-|\delta|-s}\partial_r(f^2(t,r)-1) = 0 \quad  (s=0,1,2) }
for all  $t\in [0,T_M)$ and $\delta<0$.
\end{prop}

\begin{proof}
To obtain the boundary conditions (\ref{eq4.34}), multiply
(\ref{eq4.31}) by $f^3$, (\ref{eq4.32}) by $f^2$, take the limit,
and use that $f$ is a bounded function of $r$. The fact that $f$
solves (\ref{eq4.18}), subject to these conditions, for all
$(t,r)\in[0,T_M)\times [0,\infty)$ follows from the facts that (i)
$q$ and $h$ enter (\ref{eq4.3}) which solves Ricci flow
(\ref{eq1.1}), (ii) $f$ enters (\ref{eq4.12}) which solves
Hamilton-DeTurck flow (\ref{eq1.2}), and (iii) the diffeomorphism
(\ref{eq4.8}--\ref{eq4.11}) relating these flows is valid for all
such $(t,r)$ (Proposition 4.5).
\end{proof}

%\begin{thm}
%\label{spcont} \mnote{[spcont]}
\begin{thm}\label{Thm4.7}
If there exists a constant $C_{\lambda}>0$ independent of $T_M$ such
that
%spcont1
\leqn{eq4.36}{ \sup_{0<r<\infty}\bigl(|\lambda_1(t,r)| +
|\lambda_2(t,r)|\bigr) \leq C_{\lambda} \, , }
then $T_M = \infty$.
\end{thm}

\begin{proof}
{}From Proposition (\ref{Prop4.6}), the solution $f$ of
(\ref{eq4.18}) exists up to time $T_M$. From
(\ref{eq4.27}--\ref{eq4.30}), if the sectional curvatures
$\lambda_1$ and $\lambda_2$ are bounded independent of $T_M$, then
so is $\vert {\overline{\rm Rm}} \vert^2$, and then by Theorem 3.5
we have $T_M = \infty$.
\end{proof}

\subsection{Quasi-Local Mass}

\noindent This subsection is a brief aside, not necessary for our
main results, but intended to relate our results to the motivation
discussion in the introduction.

One of the more popular quasi-local mass formulations is the
Brown-York mass. The Brown-York quasi-local mass contained within a
closed hypersurface $\Sigma$ is defined to be
\begin{equation}
\mu[\Sigma]:=\int_\Sigma(H_0-H)d\Sigma\ , \label{eq4.37}
\end{equation}
where $H$ is the mean curvature of $\Sigma$ and $H_0$ is the mean
curvature of the image of $\Sigma$ under an isometric embedding of
$\Sigma$ into flat space (assuming there is such an embedding). In
the case of a hypersphere $r=b(t)$, whose coordinate radius we will
allow to possibly change in time, we have (using (\ref{eq4.15}) and
writing $d\Omega$ to represent the canonical volume element on the
$(n-1)$-sphere)
\begin{eqnarray}
\mu(t)&=&\int_{{\mathbb S}^{n-1}}\frac{1}{b(t)} \left (1 -
\frac{1}{f(t,b(t))} \right ) b^{n-1}(t) d\Omega\nonumber \\
&=&b^{n-2}(t)\left (1 - \frac{1}{f(t,b(t))} \right ) {\rm vol}\left
( {\mathbb S}^{n-1},{\rm can}\right ) \ . \label{eq4.38}
\end{eqnarray}
Comparing to (\ref{eq4.32}), we can relate quasi-local mass to
sectional curvature by
\begin{equation}
\frac{1}{b^n(t)}\left (1 + \frac{1}{f(t,b(t))} \right )\mu(t,b(t)) =
\lambda_2(t,b(t)){\rm vol}\left ( {\mathbb S}^{n-1},{\rm can}\right
) \ . \label{eq4.39}
\end{equation}

\begin{prop}\label{Prop4.8}
The sign of the Brown-York quasi-local mass within
the hypersphere $r=b(t)$ at time $t$ is determined by the sign of
$\lambda_2(t,b(t))$, and
\begin{equation}
\lim_{t\to\infty}\lambda_2(t,b(t))=0\ \Leftrightarrow\
\lim_{t\to\infty}\mu(t,b(t))=0\ . \label{eq4.40}
\end{equation}
\end{prop}

\begin{proof}
Obvious from (\ref{eq4.39}) and (\ref{eq4.23}).
\end{proof}

Perhaps the three most interesting kinds of hyperspheres are those
of
\begin{enumerate}
\item[(i)] fixed surface area
\begin{equation}
b(t)=b_0=const>0 \ , \label{eq4.41}
\end{equation}
\item[(ii)] fixed volume contained within
\begin{equation}
\int_0^{b(t)} \int_{{\mathbb S}^{n-1}} f(t,r)r^{n-1} drd\Omega =:V_0
=const>0 {\rm \ , and} \label{eq4.42}
\end{equation}
\item[(iii)] fixed proper radius
\begin{equation}
\int_0^{b(t)}f(t,r)dr =: R_0=const>0 \ . \label{eq4.43}
\end{equation}
\end{enumerate}
In either case, it is easy to see that
\begin{equation}
0< C^-_b \le b(t) \le C^+_b \ . \label{eq4.44}
\end{equation}
where obviously $C_b^{\pm}=b_0$ for the fixed area case, while
\begin{equation}
C^{\pm}_b = \frac{nV_0}{C^{\mp}_{f^2}{\rm vol}\left ({\mathbb
S}^{n-1},{\rm can}\right )} \label{eq4.45}
\end{equation}
for the fixed volume case and
\begin{equation}
C^{\pm}_b:= \frac{R_0}{C^{\mp}_{f^2}}\label{eq4.46}
\end{equation}
for the fixed proper radius case.

\begin{rem}\label{Rem4.9}
{\rm In Subsection 5.1, we prove that $\lambda_2(t,r)\sim 1/t$ for
large $t$ and fixed $r$. Thus, for all three kinds of hyperspheres
discussed above, the quasi-local mass vanishes like $1/t$ as
$t\to\infty$.}
\end{rem}

\section{Immortality and Convergence}
\setcounter{equation}{0}

\noindent In the next two subsections we show that the sectional
curvatures $\lambda_1$ and $\lambda_2$ are bounded on $t\in[0,T_M)$.
(Equivalently, we obtain bounds on the quasi-local mass and its
radial derivative.) This permits us to invoke Theorem \ref{Thm4.7}
to conclude that the solution is immortal. In fact, we find bounds
that actually decay in time, going to zero in the limit $t\rightarrow \infty$.
This implies that the flow converges in the limit to a space with
vanishing sectional curvatures; i.e., to a flat space. In Subsection
5.3, we prove that it converges to Euclidean space ${\mathbb E}^n$.

In this section, we use $T$ to denote an arbitrary time that is less
than the maximal time of existence, i.e., $0<T<T_M$.

\subsection{The Decay of $\lambda_2$}

\noindent Short-time existence guarantees that $f^2(t,r)-1\in{\cal
O}(r^2)$ as $r\to 0$. Specifically, for all $r<r_0$ and for $0\le t
<T_M$, there is a function $C(t)$ such that
\begin{equation}
\vert w(t,r)\vert = \vert f^2(t,r)-1\vert < C(t) r^2\ .
\label{eq5.1}
\end{equation}
This follows by applying the boundedness of $f^2$ (\ref{eq4.23}) to
equation (\ref{eq4.32}) governing $\lambda_2$, which can be written
(by choosing $C(t)$ less than optimally perhaps) as
\begin{equation}
r^2\vert \lambda_2 (t,r) \vert = \left \vert\frac{1}{f^2}-1\right
\vert<C(t)r^2 \ , \label{eq5.2}
\end{equation}
To apply the continuation principle, we need to prove that $C(t)$ is
bounded in $t$. In this section we will prove more: we will show
that $C(t)$ can be taken to decay in time, converging to zero in the
limit $t\to\infty$, so that the sectional curvature $\lambda_2$
decays to zero as well.

If $w=f^2-1$ decays, then, based on the parabolic form of
(\ref{eq4.20}), one might speculate that this decay would go roughly
like $r^2/t$, or inverse ``parabolic time''. If so, then the
function $g(t,r)(f^2-1)$ should be bounded if we take $g\sim t/r^2$.
We will show below that this expectation is basically correct.

We do not take $g=t/r^2$ exactly. For small $t$, we will modify the
form $g\sim t/r^2$ so that $g$ does not vanish at $t=0$. For small
$r$, the form $g\sim t/r^2$ is problematical because we cannot
specify {\it a priori} the behaviour of $\frac{1}{r^2}(f^2-1)$ on
approach to $r=0$. This behaviour is governed by $C(t)$, the very
quantity we seek to control as the {\it outcome} of the argument, so
we cannot specify it as input. We therefore choose instead small $r$
behaviour of the form $g(t,r)\sim 1/r^m$, $m<2$, and only later do
we take $m\to 2$. For $m<2$, $g(t,r)(f^2-1)$ is very well controlled
{\it a priori} for small $r$: it goes to zero. Lastly, as
foreshadowed by (\ref{eq5.2}), we need to apply these considerations
not only with $f^2-1$ but also to $\frac{1}{f^2}-1$. The same
heuristic reasoning leads us then to consider functions of the form
$g(t,r)(\frac{1}{f^2}-1)$ with the same $g(t,r)$.

\begin{Def}\label{Def5.1}
{\em Let $f$ be defined by (\ref{eq4.13})\footnote
{wherein, of course, $q$ and $h$ arise from an asymptotically flat
Ricci flow of rotationally symmetric initial data obeying the
conditions of Theorem 1.1.}
and therefore have all the properties outlined in Section 4. For
such an $f$, define the {\em $u_m$ functions}, $m\in (0,2]$, on
$[0,T_M) \times [0,\infty)$ by
\begin{eqnarray}
u_m(t,r)&:=&\left ( \frac{1+t}{r^m+r^2}\right ) \left (
\frac{1}{f^2(t,r)}-1 \right
)\ \text{for}\ r>0\ , \label{eq5.3}\\
u_m(t,0)&:=& =\lim_{r\to 0} u_m(t,r)\ .\nonumber
\end{eqnarray}
}
\end{Def}

The $u_m$ functions have the following properties, which follow from
the flow equation (\ref{eq4.18}) for $f$, Proposition \ref{Prop4.6},
and equation (\ref{eq4.14}):
\begin{enumerate}
%\item[(i)] They are well-defined on $[0,T_M)\times [0,\infty)$.
\item[(i)] $u_m(t,0)=0$ for all $0<m < 2$ and $\lim_{r\to\infty}
u_m(t,r) =0$ for all $0<m \leq 2$.
\item[(ii)] For fixed $t$ and
$r\neq 0$, the map $m\mapsto u_m(t,r)$ is continuous at $m=2$.
\item[(iii)]
\begin{equation}
\lambda_2=-\frac{2}{1+t}u_2 \ . \label{eq5.4}
\end{equation}
\item[(iv)] The $u_m$ obey a maximum principle, as we will show
below.
\item[(v)] By direct calculation starting from
(\ref{eq4.18}), the $u_m$ obey the differential equation:
\end{enumerate}
\begin{eqnarray}
\frac{\partial u_m}{\partial t} &=& \frac{1}{f^2} \frac{\partial^2
u_m}{\partial r^2} - \frac{(r^m+r^2)}{2(1+t)} \left ( \frac{\partial
u_m}{\partial r} \right )^2-
\frac{(2r+mr^{m-1})}{(1+t)}u_m\frac{\partial u_m}{\partial r}\nonumber\\
&&+\left [ \frac{2\left ( 2r+mr^{m-1}\right )}{(r^m+r^2)f^2}
-\frac{1}{rf^2}+\frac{(n-2)}{r}
\right ] \frac{\partial u_m}{\partial r}\nonumber\\
&&-\frac{(2-m)(m+n-2)}{r^2(1+r^{2-m})}u_m\nonumber\\
&&+\frac{1}{(1+t)}\biggl \{\frac{1}{1+r^{2-m}}\left [ u_m-\left (
(4-m)(m+n-2)+m(n-2)\right )u_m^2 \right ]\nonumber\\
&& +\frac{r^{2-m}}{1+r^{2-m}}\left [ u_m-2(n-1)u_m^2\right
]\nonumber\\
&& +\frac{r^{m-2}}{1+r^{2-m}}\left [ (m-2)(m+n-2)-m\left (
\frac{m}{2}+n-2\right ) \right ]u_m^2\biggr \}\ . \label{eq5.5}
\end{eqnarray}
This PDE is the starting point for the maximum principle, which we
now derive.

\begin{prop}\label{Prop5.2}
For $u_m(t,r)$ defined by Definition \ref{Def5.1}, there is a
constant $C^+_u$ which depends only on the initial data
$a(r)=f(0,r)$ such that $u_m(t,r)\le C^+_u$ for all $t\in [0,T_M)$
and all $m\in (0,2)$.
\end{prop}

\begin{proof}
The technique will be to solve (\ref{eq4.18}) for $f$, given initial
data obeying the bounds in Corollary \ref{Cor4.3}. {}From this
initial data, we can construct initial data for $u_m$ using from
(\ref{eq5.3}) that
\begin{equation}
u_m(0,r):=\frac{1}{r^m+r^2}\left ( \frac{1}{f^2(0,r)}-1\right )=
-\frac{\lambda_2(0,r)}{1+r^{m-2}}\label{eq5.6} \ .
\end{equation}
Now by the assumed differentiability and asymptotic flatness of the
initial metric stated in Theorem \ref{Thm1.1}, the initial sectional
curvature $\lambda_2(0,r)$ is bounded. In particular, then by
(\ref{eq5.6}) $u_m(0,r)$ is bounded above on $r\in[0,\infty)$ by a
constant $C^+_u$ which depends only on the initial metric (thus on
$a(r)$ as in (\ref{eq4.4})) and so does not depend on $m$. Without
loss of generality, we choose $C^+_u\ge \frac{1}{2(n-1)}$, for
reasons that will become clear. Now it remains to be shown that
$u_m(t,r)$ is bounded above for all time $t\ge 0$ by a bound that is
dependent only on $u_m(0,r)$. Of course, the initial data
$u_m(0,r)$ will vary with $m$ (because of the denominator of
(\ref{eq5.6}); but $f(0,r)=a(r)$ and, thus, $\lambda_2(0,r)$ are of
course independent of $m$), but $C^+_u$ will always provide an
$m$-independent upper bound which will then bound the full solution.

First restrict consideration to the compact domain $D=[0,T]\times
[r_1,r_2]$, $0<r_1<r_2$, $T<T_M$, with parabolic boundary $P$ (as
defined in the proof of Proposition \ref{Prop4.2}). Now consider in
(\ref{eq5.5}) the terms that do not contain derivatives. There are
three such terms, each comprised of a function of $r$ multiplying a
factor in square brackets. One can easily check (e.g., by direct
substitution; keep in mind that $m\in(0,2)$ and $n\ge 3$) that in
(\ref{eq5.5}) each of these factors in square brackets is negative
whenever
\begin{equation}
u_m>\frac{1}{2(n-1)}\ , \label{eq5.7}
\end{equation}
so
\begin{eqnarray}
\frac{\partial u_m}{\partial t} &<& \frac{1}{f^2} \frac{\partial^2
u_m}{\partial r^2} - \frac{(r^m+r^2)}{2(1+t)} \left ( \frac{\partial
u_m}{\partial r} \right )^2-
\frac{(2r+mr^{m-1})}{(1+t)}u\frac{\partial u_m}{\partial r}\nonumber\\
&&+\left [ \frac{2\left ( 2r+mr^{m-1}\right )}{(r^m+r^2)f^2}
-\frac{1}{rf^2}+\frac{(n-2)}{r} \right ] \frac{\partial
u_m}{\partial r} \label{eq5.8}
\end{eqnarray}
then. Applying the usual maximum principle argument to this
inequality (i.e., evaluating both sides at a hypothesized local
maximum and observing that the inequality cannot then be satisfied),
we conclude that $u_m$ has no maximum greater than
$\frac{1}{2(n-1)}$ in $D\backslash P$.

By the properties of $u_m$ listed above, we have $u_m(t,r)\to 0$
both for $r\to 0$ and for $r\to \infty$. Thus, as with the proof of
Proposition \ref{Prop4.2}, if the maximum of $u_m$ is
$>\frac{1}{2(n-1)}$ (or merely positive) and lies on the parabolic
boundary with $r_1$ chosen small enough and $r_2$ large enough, it
must lie on the initial boundary. Taking the limits $r_1\to 0$ and
$r_2\to\infty$, then we see that
\begin{equation}
u_m(t,r)\le \max \left \{ \frac{1}{2(n-1)}, \sup_{r\in [0,\infty)}
\{ u_m(0,r) \} \right \} \le C^+_u \ , \label{eq5.9}
\end{equation}
for all $(t,r)\in [0,T]\times[0,\infty)$. But this holds for any
$T<T_M$, so it holds for $(t,r)\in [0,T_M)\times[0,\infty)$.
\end{proof}

\begin{cor}\label{Cor5.3}
Proposition 5.2 extends to the case $m=2$ and yields
\begin{equation}
\lambda_2(t,r)\ge -\frac{2C^+_u}{1+t} =: \frac{C^-_{\lambda_2}}{1+t}
\ . \label{eq5.10}
\end{equation}
\end{cor}

%\bigskip\noindent{\bf Proof.}
\begin{proof}
As in Proposition \ref{Prop5.2}, we solve (\ref{eq4.18}) with the
assumed initial data to find $f$, from which we construct $u_m$ for,
say, $0<m\le 2$. Fixing any $t\in[0,T_M)$ and any $r\neq 0$, the map
$m\mapsto u_m(t,r) =\frac{(1+t)}{r^2+r^m} \left (
\frac{1}{f^2(t,r)}-1\right )$ is obviously continuous at $m=2$. This
and Proposition 5.2 imply that $u_2(t,r)\le C^+_u$ for all $r>0$. By
the continuity of $r\to u_2(t,r)$, then $u_2(t,0)\le C^+_u$ as well,
for all $t\in[0,T_M)$. Now use (\ref{eq5.4}).
\end{proof}

Thus $\lambda_2$ is bounded below by a bound that tends to zero in
the limit of long times. Next we need a similarly decaying bound
from above. To get it, we work with the following class of
functions:

%\bigskip\noindent {\bf Definition 5.4.}
\begin{Def}\label{Def5.4}
{\em Let $f$ be defined by (\ref{eq4.13}) and therefore have all the
properties outlined in Section 4. For such an $f$, the {\em $v_m$
functions}, $m\in (0,2]$ are defined on $[0,T_M) \times [0,\infty)$
as
\begin{eqnarray}
v_m(t,r)&:=&\left ( \frac{1+t}{r^m+r^2}\right ) \left ( f^2(t,r)-1
\right )\ \text{for}\ r>0\ , \label{eq5.11}\\
v_m(t,0)&:=& =\lim_{r\to 0} v_m(t,r)\ .\nonumber
\end{eqnarray}
}
\end{Def}

These functions have essentially the same properties as those listed
for the $u_m$, but the relation to $\lambda_2$ is now
\begin{equation}
v_2(t,r)=\frac{1}{2}(1+t)f^2(t,r)\lambda_2(t,r)\ ,\label{eq5.12}
\end{equation}
and the $v_m$ obey the PDE (computed directly from (\ref{eq4.18})
and (\ref{eq5.11}))
\begin{eqnarray}
\frac{\partial v_m}{\partial t}&=&\frac{1}{f^2}\frac{\partial^2
v_m}{\partial r^2} -\frac{3(r^m+r^2)}{2f^4(1+t)}\left (
\frac{\partial v_m}{\partial r} \right )^2\nonumber\\
&&+\left [ \frac{2(mr^{m-1}+2r)}{(r^m+r^2)f^2} -\frac{3(mr^{m-1}+2r)}{(1+t)f^4}+\frac{n-2}{r}-\frac{1}{rf^2} \right ] \frac{\partial v_m}{\partial r}\nonumber\\
&&+\left [ 1-\frac{3(mr^{m-1}+2r)^2}{2(r^m+r^2)f^4}v_m\right ] \frac{v_m}{1+t} \nonumber\\
&&+\frac{(m-2)}{r^2} \left ( \frac{r^m}{r^m+r^2} \right )
\left ( n-2+\frac{m}{f^2} \right ) v_m\ . \label{eq5.13}
\end{eqnarray}

We must of course prove that the $v_m$ obey a maximum principle. In
fact, Proposition \ref{Prop5.2} holds with $v_m$ replacing $u_m$ and
with $m$ restricted this time to $1<m<2$. Just as with Corollary
\ref{Cor5.3}, the result can be extended to cover $m=2$. To prove
this, it will help to note that when $v_m\ge 0$, $n\ge 3$, and
$1<m<2$, then we can discard  most of the nonderivative terms in
(\ref{eq5.13}) to obtain
\begin{eqnarray}
\frac{\partial v_m}{\partial t}&\le&\frac{1}{f^2}\frac{\partial^2
v_m}{\partial r^2} -\frac{3(r^m+r^2)}{2f^4(1+t)}\left (
\frac{\partial v_m}{\partial r} \right )^2\nonumber\\
&&+\left [ \frac{2(mr^{m-1}+2r)}{(r^m+r^2)f^2} -\frac{3(mr^{m-1}+2r)}{(1+t)f^4}+\frac{n-2}{r}-\frac{1}{rf^2} \right ] \frac{\partial v_m}{\partial r}\nonumber\\
&&+\frac{v_m}{(1+t)}\left [ 1 -\frac{6v_m}{f^4} \right ]\ , \qquad v_m>0 \ .
\label{eq5.14}
\end{eqnarray}

\begin{prop}\label{Prop5.5}
There is a constant $C^+_v$ which depends only on the initial data
$f(0,r)=a(r)$ such that $v_m(t,r)<C^+_v$ for all $(t,r)\in
[0,T_M)\times [0,\infty)$ and all $m\in (1,2)$.
\end{prop}

\begin{proof}
The proof follows that of Proposition \ref{Prop5.2}. Consider first
the initial data
\begin{equation}
v_m(0,r)=\left ( \frac{1}{r^m+r^2}\right ) \left ( f^2(0,r)-1 \right
) =\left ( \frac{1}{1+r^{m-2}}\right )\frac{w(0,r)}{r^2}\ \le C^+_v
\label{eq5.15}
\end{equation}
because $\frac{|w(0,r)|}{r^2}$ is bounded, where $C^+_v$ is
independent of $m$. This time, we will choose without loss of
generality that $C^+_v\ge \frac{1}{6}(C^+_{f^2})^2$, for reasons that
will become clear below.

Again we work first on the domain $D=[0,T]\times[r_1,r_2]$,
$0<r_1<r_2$, with parabolic boundary $P$. Observe that if $v_m>
\frac{1}{6}(C^+_{f^2})^2$, the last term in (\ref{eq5.14}) will be
negative. As before, elementary arguments applied to (\ref{eq5.14})
imply that this term cannot be negative at a maximum in $D\backslash
P$, and thus such a maximum can occur only on $P$. Also as before, we
take $r_1\to 0$, $r_2\to \infty$ and since $v_m$ vanishes in both
limits, the maximum of $v_m$, if it is greater than
$\frac{1}{6}(C^+_{f^2})^2$, must occur on the initial boundary where
$t=0$. Thus we obtain for any $(t,r)\in [0,T] \times [0,\infty)$
that
\begin{equation}
v_m(t,r)\le \max \left \{\frac{1}{6}(C^+_{f^2})^2,\max_{r\in[0,\infty)}
\{v_m(0,r) \} \right \} \le C^+_v \ , \label{eq5.16}
\end{equation}
and $C^+_v$ does not depend on $m$. It also does not depend on $T$
and so the result extends to hold for all $(t,r)\in [0,T_M) \times
[0,\infty)$.
\end{proof}

\begin{rem}\label{Rem5.6}
{\rm For use in the next subsection, we observe that in virtue of
this result $u_m$ is now bounded below, as well as above, on
$(t,r)\in [0,T_M) \times [0,\infty)$ by a bound that depends only on
the initial data for $f$ and so is independent of $m$. The proof is
to observe that $u_m=-v_m/f^2\ge -C^+_v/C^-_{f^2}=:C^-_u$. We define
\begin{equation}
C_u:=\max \{ \vert C^{\pm}_u\vert \}\ , \label{eq5.17}
\end{equation}
which bounds the magnitude of $|u_m|$ and is independent of $m$.}
\end{rem}

\begin{cor}\label{Cor5.7}
Proposition \ref{Prop5.5} extends to the case $m=2$ and yields
\begin{equation}
\lambda_2(t,r)\le
\frac{2C^+_v}{C^-_{f^2}(1+t)}=:\frac{C^+_{\lambda_2}}{1+t}\ .
\label{eq5.18}
\end{equation}
\end{cor}

\begin{proof}
The extension to $m=2$ follows exactly as in Corollary \ref{Cor5.3}.
Equation (\ref{eq5.18}) follows directly from (\ref{eq5.12}).
\end{proof}

\begin{prop}\label{Prop5.8}
$\vert \lambda_2 \vert$ is bounded on $[0,T_M)\times [0,\infty)$ and
if $T_M=\infty$ then $\lambda_2$ converges uniformly to zero as
$t\to\infty$.
\end{prop}

\begin{proof}
Immediate from Corollaries \ref{Cor5.3} and \ref{Cor5.7}.
\end{proof}

\noindent In this regard, note that by Theorem \ref{Thm4.7} we {\it
can} assume $T_M=\infty$ if we can bound $\lambda_1$, which we now
proceed to do.

\subsection{The Decay of $\lambda_1$}

\noindent A lower bound and decay estimate on $\lambda_1$ is now
easy to obtain. It is quickest to work from the flow equation for
the scalar curvature, which is
\begin{eqnarray}
\frac{\partial R}{\partial t} &=& \Delta R +\xi \cdot \nabla
R + 2 R_{ij}R^{ij}\nonumber \\
&\ge&\Delta R +\xi \cdot \nabla R + \frac{2}{n} R^2\ ,
\label{eq5.19}
\end{eqnarray}
with $\xi=\xi_1dr$ given by (\ref{eq4.17}) and where we used the
elementary identity $R^{ij}R_{ij}\ge \frac{1}{n}R^2$. Inequality
(\ref{eq5.19}) gives a well-known minimum principle for $R$.
Moreover, if we define
\begin{equation}
{\tilde R}:=(1+t)R\ , \label{eq5.20}
\end{equation}
we obtain from (\ref{eq5.19}) that
\begin{equation}
\frac{\partial {\tilde R}}{\partial t} \ge \Delta {\tilde R} + \xi
\cdot \nabla {\tilde R} + \frac{1}{(1+t)} \left ( \frac{2}{n}
{\tilde R}^2 + {\tilde R} \right )\ , \label{eq5.21}
\end{equation}
which also has a minimum principle.

\begin{prop}\label{Prop5.9}
If $R$ is the scalar curvature of a Ricci flow developing from
asymptotically flat initial data on a manifold $M$ then there is a
constant $C^-_R\le 0$ such that on $[0,T_M)\times[0,\infty)\ni
(t,r)$ we have
\begin{equation}
R\ge \frac{C^-_R}{1+t} \ . \label{eq5.22}
\end{equation}
\end{prop}

For notational convenience, we give the proof for the special case
of interest, a rotationally symmetric flow on ${\mathbb R}^n$, but
the proof clearly generalizes to arbitrary asymptotically flat
flows.

\begin{proof}
First take $t\in[0,T]$, $T<T_m$. Let $B_0(a)$ be the ball of
coordinate radius $r=a$ about the origin $0\in{\mathbb R}^n$ at time
$t$. Applying elementary minimum principle arguments to
(\ref{eq5.21}), it is clear that either the minimum of ${\tilde R}$
in $[0,T]\times B_0(a)$ occurs on the parabolic boundary $P$ or
${\tilde R}\ge -\frac{n}{2}$. Now the parabolic boundary has an
initial component $t=0$ and a spatial component which is a sphere
$r=a$ for all $t>0$. By asymptotic flatness, $R\to 0$ as $a\to
\infty$ and hence ${\tilde R}\to 0$ as well. Taking this limit, we
conclude that if ${\tilde R}$ is anywhere less than $-\frac{n}{2}$,
then the minimum of ${\tilde R}$ over all $(t,x)\in [0,T]\times
{\mathbb R}^n$ exists and is realized on the initial boundary. Thus
choose $C^-_R=\min \left \{ -\frac{n}{2}, \inf_r \{ R(0,r) \} \right
\}$, which is obviously independent of $T$, so finally take $T\to
T_M$. Then ${\tilde R}\ge C^-_R \Rightarrow R\ge C^-_R (1+t)$ for
all $(t,r)\in [0,T_M)\times [0,\infty)$.
\end{proof}

\begin{cor}\label{Cor5.10}
Then $\lambda_1(t,r)$ is bounded below on $[0,T_M)\times
[0,\infty)\ni (t,r)$ by
\begin{equation}
\lambda_1(t,r)\ge \frac{1}{(1+t)} \left ( \frac{C^-_R}{2(n-1)} -
\frac{(n-2)C^+_v}{C^-_{f^2}} \right ) =: \frac{C^-_{\lambda_1}}{1+t}
\ . \label{eq5.23}
\end{equation}
\end{cor}

\begin{proof}
This follows from the formula
\begin{equation}
R=2(n-1)\lambda_1+(n-1)(n-2)\lambda_2 \label{eq5.24}
\end{equation}
for the scalar curvature in terms of the sectional curvatures,
equation (\ref{eq5.22}), and the upper bound (\ref{eq5.18}) on
$\lambda_2$.
\end{proof}

Now we turn attention to finding an upper bound and decay estimate.
We have to work harder than we did for the lower bound, but we can
apply essentially the same strategy as we used to prove boundedness
and convergence of $\lambda_1$. Once again, the main issue will be
control of $\lambda_1$ at $r=0$, and we will be forced to work with
a sequence of functions with known behaviour at $r=0$. This time, we
have found that a choice well-suited to our purpose is given by

\begin{Def}\label{Def5.11}
{\em Let $f$ be defined by (\ref{eq4.13}) and therefore have all the
properties outlined in Section 4. For such an $f$, define the {\em
$y_m$ functions}, $m\in \left (1,2\right ]$, on $[0,T_M) \times
[0,\infty)$ by
\begin{eqnarray}
y_m(t,r)&:=&\left ( \frac{1+t}{1+r^{2-m}}\right ) \left \{ r
\frac{\partial}{\partial r}\left [ \frac{1}{r^m} \left (
\frac{1}{f}-1\right ) \right ] \right \},\quad r>0, \label{eq5.25}\\
y_m(t,0)&:=&\lim_{r\to 0} y_m(t,r)\ . \nonumber
\end{eqnarray}
}
\end{Def}

We can extract $\lambda_1$ from the relation
\begin{equation}
\frac{y_m}{1+t}=\frac{r^2f}{(r^m+r^2)}\left (
\frac{m}{(1+f)}\lambda_2 - \lambda_1 \right ) \ . \label{eq5.26}
\end{equation}

Notice that $y_m(t,r)\to 0$ as $r\to 0$ whenever $m<2$. Calculating
from (\ref{eq4.18}), we find that $y_m(t,r)$ obeys
\begin{eqnarray}
\frac{\partial y_m}{\partial t} &=& \frac{1}{f^2}\frac{\partial^2
y_m}{\partial r^2}+\frac{1}{r}\alpha_m\frac{\partial y_m}{\partial
r} \nonumber\\
&&+\frac{1}{r^2}\biggl \{ \left [ \frac{2}{f}\left ( 2(m-1)r^m
+mr^2\right ) y_m +1\right ] \frac{y_m}{1+t} \nonumber\\
&&+\beta_m y_m+(1+t)\gamma_m\biggr \} \ , \label{eq5.27}
\end{eqnarray}
where some of the coefficients have rather lengthy expressions so we
have introduced the abbreviations
\begin{eqnarray}
\alpha_m&:=&
\frac{2(r^{m}+r^2)}{f}\frac{y_m}{(1+t)}+\frac{4m-3}{f^2}
-\frac{2m}{f}+n-2\nonumber\\
&&-\frac{2(m-2)r^{2-m}}{f^2(1+r^{2-m})} \ ,\label{eq5.28} \\
\beta_m&:=&\frac{7m^2-14m+4}{f^2} -\frac{m(6m-8)}{f} + (n-2)\left
( m-1-\frac{3}{f^2} \right ) \nonumber\\
&&+\frac{(m-2)r^{2-m}}{1+r^{2-m}}\left [ -\frac{(3m-2)}{f^2}
+ \frac{2m}{f}-(n-2)\right ]\ , \label{eq5.29} \\
\gamma_m&:=&\frac{1}{(r^m+r^2)} \left ( \frac{1}{f} -1\right )
\biggl \{
\frac{2m(m-1)(m-2)}{f^2} +\frac{2m^2(2-m)}{f}\nonumber\\
&&+(n-2)\left [ -m+\frac{m+2}{f}+\frac{2(1-m)}{f^2}\right ] \biggr
\} \ . \label{eq5.30}
\end{eqnarray}

We now claim that the $y_m(t,r)$ are bounded below on $[0,T_M)\times
[0,\infty)$ by a constant that is independent of $m$. Proceeding in
our now usual fashion, let $T$ be such that $0<T<T_M$ and define
$D:=[0,T]\times[0,\infty)$. As usual, because $y_m$ tends to zero
for $r\to\infty$ and for $r\to 0$, either zero is the lower bound or
\begin{equation}
\inf_D y_m=:Y=y_m(t_0,r_0)<0 \label{eq5.31}
\end{equation}
for some $t_0$ and some $r_0>0$. In the latter case, either $t_0=0$
and therefore the minimum depends only on initial data $a(r)=f(0,r)$
and not on $m$ or $T$, or it occurs at some $t_0\in (0,T]$ and then
the minimum obeys a quadratic inequality which we now state:

\begin{lem}\label{Lem5.12}
Let $y_m$ be defined on $[0,T]\times [0,\infty)$, $T<T_M$, by
Definition \ref{Def5.11}. For $m<2$, if $y_m$ has a negative infimum
$Y<0$, then this infimum is realized as a minimum at some
$(t_0,r_0)$ where $r_0>0$ and either $t_0=0$ or
\begin{equation}
\frac{2(m-1)}{f(t_0,r_0)}Y^2+\left ( \frac{1+t_0}{r_0^m+r_0^2}
\right ) \left [ \left ( 1+\beta_m(t_0,r_0)\right
)Y+(1+t_0)\gamma(t_0,r_0) \right ] <0 \ . \label{eq5.32}
\end{equation}
\end{lem}

\begin{proof}
As discussed immediately above, a negative infimum must be realized
at some $(t_0,r_0$ where $r_0>0$. Then it follows by applying
standard minimum principle arguments to equation (\ref{eq5.27}) that
either the minimum occurs at $t_0=0$ or the nonderivative terms in
(\ref{eq5.27}) are governed by the inequality
\begin{eqnarray}
0&>&\frac{2}{f(t_0,r_0)}\left ( 2(m-1)r_0^m +mr_0^2\right ) Y^2 +Y
\nonumber\\
&&+(1+t_0)\left [ \beta_m(t_0,r_0) Y+(1+t_0)\gamma_m(t_0,r_0)\right
]\ . \label{eq5.33}
\end{eqnarray}
But in the first term on the right-hand side, use that $\left (
2(m-1)r_0^m +mr_0^2\right ) Y^2 > 2(m-1)\left ( r_0^m+r_0^2\right
)Y^2$ for $m<2$ to replace the former by the latter. Replace the
second term (the singleton $Y$) by $(1+t_0)Y < Y$. These
replacements preserve the inequality. Divide by $r_0^m +r_0^2$ to
complete the proof.
\end{proof}

Now further restrict $m$ to some range of form $1<\kappa\le m<2$, so
that the coefficient of $Y^2$ in (\ref{eq5.32}) is not arbitrarily
small; for definiteness $\kappa=\frac{3}{2}\le m <2$ will do nicely.
Then since the criterion (\ref{eq5.32}) is quadratic in $Y$ with
positive coefficient of $Y^2$, it will be violated for $Y$
sufficiently negative. Thus $Y$ cannot be arbitrarily negative,
giving a bound on $y_m$ expressed in terms of the coefficients in
(\ref{eq5.32}). It remains therefore to manipulate these
coefficients to produce a bound that is manifestly independent of
$m$ and $T$. The proof is an exercise in elementary manipulation,
but we will give the main points.

\begin{prop}\label{Prop5.13}
Let $\frac{3}{2}\le m<2$. Then for each $m$, the $y_m$ are bounded
below on $[0,T_M) \times [0,\infty)$ by an $m$-independent constant.
\end{prop}

\begin{proof}
As usual, we work on $t\in[0,T]$ with $T<T_M$ to obtain a bound
which does not depend on $m$ or $T$ and then take $T\to T_M$ when
we're done.

If the lower bound is zero, which occurs at $r_0=0$ and as
$r_0\to\infty$, then obviously it is independent of $m$ and $T$, so
assume that the lower bound is negative. Then it is realized as a
minimum at some $(t_0,r_0)\in[0,T]\times[0,\infty)$. If $t_0=0$, the
lower bound is given by the initial data, so again it is clearly
$m$- and $T$-independent. Therefore, assume $t_0>0$. Then the
criterion (\ref{eq5.32}) applies.

In this last case, we start with (\ref{eq5.32}) and seek to
re-express, where possible, factors of the form
$\frac{1+t_0}{r^m_0+r^2_0}$ in terms of the bounded quantity
$u_m=\left ( \frac{1+t}{r^m+r^2}\right ) \left ( \frac{1}{f^2}-1
\right )$. The boundedness of this quantity is described in Remark
\ref{Rem5.6}; since $f$ is also bounded, we can also make use of
equivalent form $\frac{f}{1+f}u_m=\left ( \frac{1+t}{r^m+r^2}\right
) \left ( \frac{1}{f}-1 \right )$. For example, the term in
(\ref{eq5.32}) that is constant in $Y$ can be written as
(understanding all quantities to be evaluated at $(t_0,r_0)$)
\begin{eqnarray}
\frac{(1+t_0)^2}{r^m_0+r^2_0}\gamma_m&=& \left ( \frac{fu_m}{1+f}
\right )^2 \biggl [ 2m(m-2)\left ( \frac{m-1}{f}-1 \right )
\nonumber\\
&&\qquad-(n-2)\left ( \frac{2(m-1)}{f}+m-4\right ) \biggr ]
\nonumber \\
&&-2(m-2)(m+n-2)\left ( \frac{1+t_0}{r^m_0+r^2_0}\right )
\frac{fu_m}{1+f}\ . \label{eq5.34}
\end{eqnarray}
We can minimize the term proportional to $u_m^2$ over
$\frac{3}{2}\le m \le 2$. In the second term, note that the
coefficient $-2(m-2)(m+n-2)$ is positive for $\frac{3}{2}\le m < 2$.
Therefore we write $-2(m-2)(m+n-2) \frac{fu_m}{1+f} \ge \left (
\frac{1}{2}-n\right )\frac{f}{1+f} |u_m(t_0,r_0)| \ge \left (
\frac{1}{2}-n\right )\frac{f}{1+f}C_u \ge \left (
\frac{1}{2}-n\right )C_u$, using (\ref{eq5.17}). This yields
\begin{eqnarray}
\frac{(1+t_0)^2}{r^m_0+r^2_0}\gamma_m&\ge& \left ( \frac{f u_m}{1+f}
\right )^2\left [ -\frac{4}{3\sqrt{3}f}-2(n-2)\left ( 1-\frac{1}{f}
\right ) \right ]\nonumber\\
&&+\left ( \frac{1}{2}-n\right ) C_u \left (
\frac{1+t_0}{r^m_0+r^2_0} \right )\nonumber\\
&\ge&-k_1+\left ( \frac{1}{2}-n\right )C_u \left (
\frac{1+t_0}{r^m_0+r^2_0} \right )\ , \label{eq5.35}
\end{eqnarray}
where $k_1$ is a (positive) constant independent of $m$, $T$, and
$Y$.\footnote
{For example, $k_1=C_u^2\left [ \frac{1}{C^-_f} +2(n-2)\right ]$
would do fine, where we write $C^-_f:=\sqrt{C^-_{f^2}}$.}
The second term still contains an unwanted factor of
$\frac{1+t_0}{r^m_0+r^2_0}$ with negative coefficient, but for $Y$
sufficiently negative we will be able to dominate this term with
positive contributions coming from the part of the criterion
(\ref{eq5.32}) that is linear in $Y$.

To examine the linear term, start from the expression
\begin{eqnarray}
\left (\frac{1+t_0}{r_0^m+r_0^2}\right ) \left ( 1+\beta_m \right )
Y &=&\left (\frac{1+t_0}{r_0^m+r_0^2}\right )\biggl \{ 1+
\frac{7m^2-14m+4}{f^2} -\frac{m(6m-8)}{f}\nonumber\\
&& +(n-2)\left ( m-1-\frac{3}{f^2} \right ) \label{eq5.36}\\
&&+\frac{(2-m)r^{2-m}}{1+r^{2-m}} \left [ \frac{3m-2}{f^2}
-\frac{2m}{f}+n-2 \right ] \biggr \} Y\ . \nonumber
\end{eqnarray}
The terms in the last line simplify since we can use that $Y<0$,
$\frac{3}{2}\ge m < 2$, and $n\ge 3$ to write
\begin{eqnarray}
\frac{(2-m)r^{2-m}}{1+r^{2-m}} \left [ \frac{3m-2}{f^2}
-\frac{2m}{f}+n-2 \right ] Y &>& \frac{(2-m)r^{2-m}}{1+r^{2-m}}
\left [ \frac{3m-2}{f^2} +n-2 \right ] Y\nonumber \\
&>& (2-m)\left [ \frac{3m-2}{f^2} +n-2 \right ] Y\ . \label{eq5.37}
\end{eqnarray}
Now we can combine this result with (\ref{eq5.36}) and again absorb
the factor of $\frac{1+t_0}{r^m_0+r^2_0}$, wherever possible, using
$u_m$. We get
\begin{eqnarray}
\left (\frac{1+t_0}{r_0^m+r_0^2} \right ) \left ( 1+\beta_m \right
) Y &>& \frac{f u_m}{(1+f)} \biggl [ \frac{4m^2-6m-3(n-2)}{f}\nonumber\\
&&\qquad -2m^2 +2m -3(n-2) \biggr ]Y \nonumber\\
&&-\left ( 2m^2-2m +2n -5 \right ) \left ( \frac{1+t_0}{r^m_0+r^2_0}
\right ) Y\nonumber\\
&\ge&\frac{fu_m}{(1+f)} \biggl [ \frac{4m^2-6m-3(n-2)}{f}\nonumber\\
&&\qquad -2m^2 +2m -3(n-2) \biggr ] Y \nonumber\\
&&-\left ( 2n-\frac{7}{2}\right ) \left ( \frac{1+t_0}{r^m_0+r^2_0}
\right ) Y \ , \label{eq5.38}
\end{eqnarray}
where in the last line minimized over $\frac{3}{2}\le m < 2$. It is
again evident that this is the sum of a bounded term and a term
involving $\frac{1+t_0}{r^m_0+r^2_0}$. Both these terms are linear
in $Y$. That is,
\begin{equation}
\left (\frac{1+t_0}{r_0^m+r_0^2} \right ) \left ( 1+\beta_m \right )
Y\ge k_2Y-\left ( 2n-\frac{7}{2}\right ) \left (
\frac{1+t_0}{r^m_0+r^2_0} \right ) Y \ , \label{eq5.39}
\end{equation}
where $k_2$ is a constant independent of $m$, $T$, and $Y$.\footnote
{For example, from elementary considerations applied to
(\ref{eq5.38}) we obtain that $k_2=8C_u$ is a suitable bound.}

Inserting (\ref{eq5.35}) and (\ref{eq5.39}) into the criterion
(\ref{eq5.32}) and using that $\frac{2(m-1)}{f(t_0,r_0)}Y^2\ge
\frac{3}{f(t_0,r_0)}Y^2$ for $\frac{3}{2}\le m <2$, we obtain the
following necessary condition for $Y<0$ to be the minimum of
$y_m(t_0,r_0)$ at some $t_0>0$:
\begin{eqnarray}
0&\ge&\frac{3}{f(t_0,r_0)}Y^2+k_2 Y -k_1\nonumber\\
&&+\left [\left ( \frac{1}{2}-n \right ) C_u - \left (
2n-\frac{7}{2} \right )Y \right ] \left ( \frac{1+t_0}{r^m_0+r^2_0}
\right ) \ . \label{eq5.40}
\end{eqnarray}
Then a necessary condition for $Y<\frac{1-2n}{4n-7}C_u$ to be the
minimum of $y_m(t_0,r_0)$ at some $t_0>0$ is
\begin{equation}
0>\frac{3}{f(t_0,r_0)}Y^2+k_2 Y -k_1\ , \label{eq5.41}
\end{equation}
which is clearly violated whenever
\begin{equation}
Y<C_Y:=\min \left \{ \left ( \frac{1-2n}{4n-7}\right ) C_u,
-\frac{C^+_f}{6}\left [ k_2 + \sqrt{k_2^2+\frac{12}{C^-_f} k_1}
\right ] \right \} \ , \label{eq5.42}
\end{equation}
where we use the short-hand $C^{\pm}_f:=\sqrt{C^{\pm}_{f^2}}$. We
conclude that
\begin{equation}
y_m\ge C_y^-:=\min\left \{ C_Y,\inf_r \{ y_m(0,r) \} \right \}
\label{eq5.43}
\end{equation}
on $[0,T]\times [0,\infty)$ and since these bounds do not depend on
$T$, taking $T\to T_M$ we see that they hold as well on
$[0,T_M)\times [0,\infty)$
\end{proof}

\begin{cor}\label{Cor5.14}
There is a constant $C^+_{\lambda_1}$ such that
\begin{equation}
\lambda_1(t,r)\le \frac{C^+_{\lambda_1}}{1+t} \label{eq5.44}
\end{equation}
on $[0,T_M)\times [0,\infty)$.
\end{cor}

\begin{proof}
First we prove that $y_2$ is bounded below by $C_y^-$. As with
Corollaries \ref{Cor5.3} and \ref{Cor5.7}, the map $m\mapsto
y_m(t,r)$, with fixed $t$ and fixed $r>0$, is continuous, so the
bound (\ref{eq5.43}) applies to $y_2(t,r)$ except possibly at $r=0$.
Then the continuity of $y_2$ at $r=0$ implies that the bound holds
there as well.

Next, the $m=2$ case of (\ref{eq5.26}) yields
\begin{equation}
\lambda_1=\frac{2}{1+f}\lambda_2-\frac{y_2}{1+t}\ . \label{eq5.45}
\end{equation}
Using (\ref{eq5.18}) and the facts that $C^-_y\le 0$ and $C^+_v\ge
0$, we can write this as
\begin{eqnarray}
\lambda_1(t,r)&\le& \left (\frac{2}{1+f}\right
)\frac{2C^+_v}{(1+t)C^-+{f^2}} -\frac{2C_y^-}{(1+t)f}\nonumber \\
&\le& \frac{1}{1+t} \left [ \frac{4C^+_v}{C^-_{f^2}}
-\frac{2C^-_y}{C^-_f} \right ] \label{eq5.46}\ ,
\end{eqnarray}
where we've used that $0<C^-_{f^2}\le f^2$ and $C^-_f:=
\sqrt{C^-_{f^2}}$. Now let $C^+_{\lambda_1}$ equal the quantity in
square brackets in the last line.
\end{proof}

We can now prove the main theorem.

\subsection{Proof of Theorem 1.1}

%\noindent{\it Proof of Theorem \ref{Thm1.1}.}

\noindent {\it Proof of Statement (i).} By Corollaries \ref{Cor5.3},
\ref{Cor5.7}, \ref{Cor5.10}, and \ref{Cor5.14}, the sectional
curvatures in $[0,T_M)\times [0,\infty)$ are bounded above and below
by bounds of the form
\begin{equation}
\vert \lambda_{1,2}\vert \le \frac{\vert
C^{\pm}_{\lambda_{1,2}}\vert }{1+t} \le \vert
C^{\pm}_{\lambda_{1,2}}\vert\ .\label{eq5.47}
\end{equation}
Thus, by Theorem \ref{Thm4.7}, we can take $T_M=\infty$ and can
conclude that there is a constant $C_0$ such that
%curvA
\leqn{eq5.48}{ \sup_{x\in\Rbb^n}|{\rm
\overline{Rm}}(x,t)|_{\bar{g}(t,x)} \leq \frac{C_0}{1+t} \quad
\forall \; t\geq 0. }
This proves the existence for all $t\in[0,\infty)$ of the solution
developing from the initial condition Statement (i) of the theorem
and also the $\ell =0$ estimate of (iii).
\medskip

\noindent {\it Proof of Statement (ii).} This is immediate from
Theorem \ref{LocA}.
\medskip

\noindent {\it Proof of Statement (iii).} Follows directly
from \eqref{eq5.48} and Theorem 7.1 of \cite{Ham95}.
\medskip

\noindent{\it Proof of Statement (iv).} This follows from the
Compactness Theorem 1.2 of \cite{Hamilton3} and statement (iii),
provided the injectivity radius at the origin is $>\delta\ge 0$ for
some $\delta$ independent of $t$. Since the metric is uniformly
equivalent to the Euclidean metric and the sectional curvatures are
uniformly bounded in time, this follows immediately from, for
example, the Cheeger-Gromov-Taylor injectivity radius estimate
(Theorem 4.7 of \cite{CGT}).\footnote
{Even more simply, since the constant-$r$ surfaces are convex
througout the flow, there are no closed geodesics. Then it follows
from the sectional curvature bound (\ref{eq5.48}) that ${\rm inj\
}({\mathbb R^n},g(t))\ge \frac{\pi \sqrt{1+t}}{\sqrt{C_0}}$. Since
this gives a global bound on the injectivity radius, less powerful
convergence theorems (e.g Theorem 7.1.3 of \cite{Topping}) suffice
to finish the proof.}
%

%%EW: Removed:
\begin{comment}
\noindent{\it Proof of Statement (iv).} Fix any $t$ and any
$r=const>0$ hypersphere in $({\mathbb R^n},g(t))$. By rotational
symmetry, every point on this hypersurface is umbilic: all of the
principal curvatures are equal. By the absence of minimal
hyperspheres (Corollary \ref{Cor4.4}), the mean curvature is
positive and thus so is each principal curvature. Now assume a
closed geodesic exists in $({\mathbb R^n},g(t))$. The there would be
a point where the $r$-coordinate along the geodesic were a maximum,
and there the geodesic would be tangent to a constant-$r$
hypersphere. This is impossible, since one of the principal
curvatures would have to be $\le 0$ at the point of tangency. Thus,
there are no closed geodesics in $({\mathbb R^n},g(t))$ for any
$t\ge 0$.

For a manifold with no closed geodesics and with sectional
curvatures bounded above by a constant $k>0$, the injectivity radius
is bounded below by $\pi/\sqrt{k}$. In fact, in the present case, we
have from (\ref{eq5.48}) that
\begin{equation}
{\rm inj\ }({\mathbb R^n},g(t))\ge \frac{\pi \sqrt{1+t}}{\sqrt{C_0}}
\ . \label{eq5.54}
\end{equation}
{}From Statement (iii) and (\ref{eq5.54}), for {\it any} increasing,
divergent sequence of times $t_i$ starting from some $t_0>0$ and for
any convergent sequence of points $p_i\in {\mathbb R}^n$ (say choose
all $p_i$ to be the origin), conditions 7.1.1 and 7.1.2 of
\cite{Topping} apply. Then convergence follows from Theorem 7.1.3 of
\cite{Topping}.
\medskip
\end{comment}

\noindent{\it Proof of Statement (v).} Immediate from Remark
\ref{mass}. \hfill \qed

%%%% \end{document}

\appendix

%############ winq.tex ##########################33

\sect{winq}{Weighted calculus inequalities}
\setcounter{equation}{0}

In this and the following sections $C$ will be used to denote a
constant that may change value from line to line but whose exact
value will not be needed.

The next lemma is a weighted version of H\"{o}lder's inequality
and can be proved easily from the definition of the weighted $L^p$
norms and H\"{o}lder's inequality.
\begin{lem} \label{Holder} \mnote{[Holder]}
If $u\in L^q_{\delta_1}$, $v\in L^r_{\delta_2}$, $\delta = \delta_1+\delta_2$,
$1\leq p,q,r \leq \infty$, and $1/p=1/q+1/r$ then
\eqn{Holder.1}{ \norm{uv}_{L^p_\delta} \leq
\norm{u}_{L^q_{\delta_1}} \norm{v}_{L^r_{\delta_2}}\, . }
\end{lem}
It has been shown in \cite{Bart86} that the Sobolev
inequalities
extend to the weighted spaces $W^{k,p}_{\delta}$. The
one of most interest to us is:
\begin{lem} \label{SobA} \mnote{[SobA]}
%$\;\;\;$ \\
%\begin{itemize}
%\item[(ii)] If $n-kp = 0 $ and $p\leq q < \infty$ then
%$\norm{u}_{L^q_\delta} \leq C\norm{u}_{W^{k,p}_\delta}$
%for all $u\in W^{k,p}_\delta$.
%\item[(i)]
If $n- kp < 0$ then
$\norm{u}_{L^\infty_{\delta}} \leq C\norm{u}_{W^{k,p}_\delta}$
for all $u\in W^{k,p}_\delta$. Moreover $u\in C^{0}_\delta$ and
$u(x)=\text{o}(|x|^\delta)$ as $|x|\rightarrow \infty$.
%\item[(ii)]
%If $0 < \alpha \leq k-n/p \leq 1$ then
%$\norm{u}_{C^{0,\alpha}_{\delta}} \leq C\norm{u}_{W^{k,p}_\delta}$
%for all $u\in W^{k,p}_\delta$.
%\item[(iii)] If $n-kp > 0$ and $p\leq q\leq np/(n-kp)$, then
%$\norm{u}_{L^q_\delta} \leq C\norm{u}_{W^{k,p}_\delta}$.
%for all $u\in W^{k,p}_\delta$.
%\end{itemize}
\end{lem} Also of use is the weighted
multiplication lemma which follows from the weighted Sobolev
and H\"{o}lder inequalities. \begin{lem} \label{SobB} \mnote{[SobB]} If there
exists a multiplication $V_{1}\times V_{2}\rightarrow V_3 : (u,v)\mapsto
u\cdot v$ then for $1\leq p < \infty$ the corresponding
multiplication
\eqn{SobB1}{ W^{k_{1},p}_{\delta_{1}}\times W^{k_{2},p}_{\delta_2}
\rightarrow W^{k_{3},p}_{\delta_{3}} \; : \: (u,v) \mapsto u\cdot v
}
is bilinear and continuous if $k_{1},k_{2}\geq k_{3}$, $k_{3} <
k_{1}+k_{2} - n/p$, and $\delta_{1}+\delta_{2}\leq \delta_{3}$.
\end{lem}

We now introduce some notation. Let $B_R$ denote the ball of radius
$R$, and $a_R$ and $A_R$ denote the annuli $B_{2R}\setminus B_{R}$
and $B_{4R}\setminus B_{R}$, respectively. If we let
\eqn{uscal}{ u_{R}(x):= u(Rx)\, ,} then the identity \leqn{dscale}{
(D^{I}u)_{R} = R^{-|I|}D^{I}u_{R}}
and a simple change of variables argument show that for fixed
$\Lambda \geq 0$ there exists a $C>0$ independent of $u\in
W^{k,p}_{\delta}$ and $R \geq 1$ such that
\leqn{scale}{ C^{-1}R^{-\delta}\norm{u_{R}}_{W^{k,p}(a_{\Lambda})}
\leq \norm{u}_{W^{k,p}_{\delta}(a_{\Lambda R})} \leq C
R^{-\delta}\norm{u_{R}}_{W^{k,p}(a_{\Lambda})}\, \quad (1\leq p <
\infty).}
As discussed in \cite{Bart86} (see also \cite{ChrDel03}), this
inequality is the key to proving the weighted Sobolev inequalities
and weighted elliptic estimates by making it possible to turn
local estimates on $a_{\Lambda}$ and $B_{\Lambda}$ into global
ones for the weighted spaces $W^{k,p}_{\delta}$. The localization
of the scaling inequality is best seen by using it to estimate the
norm on $W^{k,p}$ in terms of local estimates and scaling. The
identity $\norm{u}^{p}_{W^{k,p}_{\delta}} =
\norm{u}^{p}_{W^{k,p}_{\delta}(B_{\Lambda})} +\sum_{j=1}^{\infty}
\norm{u}^{p}_{W^{k,p}_{\delta}(a_{2^{j-1}\Lambda})}$ and the
scaling inequality \eqref{scale} implies that there exists a $C>0$
independent of $u\in W^{k,p}_\delta$ such that
\leqn{scaleA}{C^{-1}\norm{u}^{p}_{W^{p,k}_{\delta}} \leq
\norm{u}^{p}_{W^{k,p}(B_{\Lambda})}+\sum_{j=1}^{\infty}2^{-p\delta(j-1)}
\norm{S_{j}u}^{p}_{W^{k,p}(a_{\Lambda})} \leq
C\norm{u}^{p}_{W^{k,p}_{\delta}} \quad (1\leq p < \infty)}
where $S_j$ is the scaling operator
\leqn{sop}{ S_j(u)(x):= u(2^{j-1}x). }
For $p=\infty$, it follows easily from the definition of the
weighted norm that
\leqn{scaleB}{ C^{-1}\norm{u}_{W^{k,\infty}_{\delta}} \leq
\sup\{\norm{u}_{W^{k,\infty}(B_{\Lambda})}, 2^{-\delta(
j-1)}\norm{S_j u}_{W^{k,\infty}(a_{\Lambda})}\} \leq
C\norm{u}_{W^{k,\infty}_{\delta}}}
for all $u\in W^{k,p}_\infty$.

Following \cite{Max04}, we will use a partition of
unity to express \eqref{scaleA}
in a different but equivalent form. Let $\{\phi_j\}_{j=0}^{\infty}
\subset \Co(\Rbb^n)$
be a smooth partition of unity that satisfies
\alin{punity}{
&\text{supp}\,\,\phi_0 \subset B_2\; \quad \text{supp}\, \,\phi_j \subset
A_{2^{j-1}} \quad j\geq 1\, , \\
&\phi_j(x) = \phi_1(2^{1-j}) \, .
}
Then it is not difficult to show that \eqref{scaleA} implies that
\leqn{scaleC}{C^{-1}\norm{u}^{p}_{W^{p,k}_{\delta}} \leq
\norm{\phi_0u}^{p}_{W^{k,p}}+\sum_{j=1}^{\infty}2^{-p\delta(j-1)}
\norm{S_{j}(\phi_j u)}^{p}_{W^{k,p}} \leq
C\norm{u}^{p}_{W^{k,p}_{\delta}} \quad (1\leq p < \infty)}
for all $u\in W^{k,p}_\delta$. Although in this article we will
assume that $k$ is a non-negative integer, the advantage of
\eqref{scaleC} over \eqref{scaleA} is that it can be used to define
the weighted spaces $W^{k,p}_\delta$ for non-integral $k$. See
\cite{Max04} for details.

\subsect{Moser}{Moser inequalities on the $H^{k}_{\delta}$ spaces}

We will now use the scaling technique to extend the
Moser calculus inequalities on the ordinary Sobolev spaces $H^{k}$
to the weighted ones $H^{k}_{\delta}$. For related results
on the $H^{k}_\delta$ spaces see \cite{Max04}. Also, see
\cite{ChrLen} for similar results on  a different class of
weighted Sobolev spaces.
The proof of the following lemma exemplifies the use of scaling to
get global inequalities from local ones.

\begin{lem} \label{MoserA} \mnote{[MoserA]} Suppose that $F\in
C^{\ell}_{b}(V,\Rbb)$, $F(0)=0$, and $1\leq k \leq \ell$. Then
there exists a  $C > 0$ such that
\eqn{MoserA1}{ \norm{F(u)}_{H^{k}_{\delta} }\leq C
(1+\norm{u}^{\ell-1}_{L^{\infty}})\norm{u}_{H^{k}_{\delta}}}
for all $u \in H^{k}_{\delta}\cap L^{\infty}$.\end{lem}
\begin{proof}
We recall the standard Moser inequality
(see \cite{TayIII}, Proposition 3.9)
\leqn{MoserA4}{ \norm{F(v)}_{H^{k}} \leq
C(1+\norm{v}^{\ell-1}_{L^{\infty}})\norm{v}_{H^{k}}\, \text{for all
$v\in L^\infty\cap H^k$} }
Assume $j\geq 1$. Then from the definition of the partition of unity
$\{\phi_j\}$, we have
\eqn{MoserA5}{ \phi_j F(u) = \phi_j F(\sum_{k=-1}^{+1}\phi_{j+k}u)
\, . } So using $S_j(F(u))=F(S_j(u))$, we see that \eqn{MoserA6}{
S_j(\phi_j F(u)) = \phi_1 F(S_2\circ S_{j-1}(\phi_{j-1}u) +
S_{j}(\phi_{j}u) + S_0\circ S_{j+1}(\phi_{j+1}u) )\, , } and hence
\leqn{MoserA7}{ \norm{S_j(\phi_j F(u))}^2_{H^k} \leq
C\norm{\phi_1}_{W^{k,\infty}}
(1+\norm{u}_{L^{\infty}}^{\ell-1})^2\sum_{k=-1}^{1}
\norm{S_{j-k}(\phi_{j-k}u)}^2_{H^k} }
by \eqref{MoserA4} and
\eqn{MoserA8}{\norm{S_2\circ S_{j-1}(\phi_{j-1}u) + S_{j}(\phi_{j}u)
+ S_0\circ S_{j+1}(\phi_{j+1}u)}_{L^{\infty}}
=\norm{\sum_{k=-1}^{+1}\phi_{j+k}u}_{L^\infty}\leq
\norm{u}_{L^\infty}. } Thus \eqn{MoserA9}{ \sum_{j-1}^{\infty}
2^{-2\delta(j-1)}\norm{S_j(\phi_j F(u))}^2_{H^k} \leq
C\norm{u}^2_{H^k_\delta} }
for some constant $C>0$ independent of $u\in W^{k,p}$ by
\eqref{scaleC} and \eqref{MoserA7}. Similar calculations give
\eqn{MoserA10}{ \norm{\phi_0F(u)}_{H^k}^2\leq
C\norm{u}^{2}_{H^k_\delta}. }
These two inequalities along with \eqref{scaleC} show that
$\norm{F(u)}_{H^k_{\delta}} \leq C(1+\norm{u}_{L^\infty}^{\ell-1})
\norm{u}_{H^k_{\delta}}$ for all $u\in L^{\infty}\cap H^{k}_\delta$.
\end{proof}

The other Moser inequalities that will be needed  are the weighted
version of the product and commutator estimate.
\begin{lem} \label{MoserC} \mnote{[MoserC]}
For all $|I|\leq k$, $u\in H^{k}_{\delta_{1}}\cap
L^{\infty}_{\delta_{2}}$ and $v\in H^{k}_{\delta_{3}}\cap
L^{\infty}_{\delta_{4}}$ with $\delta = \delta_{1}+\delta_{4} =
\delta_{2}+\delta_{3}$ there exist a  $C
> 0$ such that
\gath{MoserC1}{ \norm{uv}_{H^k_{\delta-1}} \leq
C\norm{u}_{H^k_{\delta_1}}\norm{v}_{L^\infty_{\delta_4}}+
C\norm{u}_{L^\infty_{\delta_2}}\norm{v}_{H^k_{\delta_3}}
\intertext{and} \norm{[D^{I},u]v}_{L^{2}_{\delta-|I|}}\leq
C\norm{
u}_{H^{k}_{\delta_{1}}}\norm{v}_{L^{\infty}_{\delta_{4}}} +
C\norm{u}_{L^{\infty}_{\delta_{2}}}\norm{v}_{H^{k}_{\delta_{3}}} \, . }
\end{lem}
\begin{proof} As in Lemma \ref{MoserA}, the proof follows
from scaling and the standard estimates (see \cite{TayIII},
Proposition 3.7) $\norm{uv}_{H^k}\leq
C\norm{u}_{H^k}\norm{v}_{L^\infty} +
C\norm{u}_{L^\infty}\norm{v}_{H^k}$ and
$\norm{[D^{I},u]v}_{L^{2}}\leq
C\norm{\nabla u}_{H^{k-1}}\norm{v}_{L^{\infty}} +
C\norm{\nabla u}_{L^{\infty}}\norm{v}_{H^{k-1}}$.
\end{proof}

In addition to the Moser inequalities, we also need to know when
the map $u \mapsto F(u)$ is locally Lipschitz on $H^{k}_{\delta}$.
\begin{lem}\label{Lip} \mnote{Lip} Suppose $F \in
C^{\ell}_{b}(V,\Rbb)$, $F(0)=0$, $\delta \leq 0$, and $k \leq \ell
$, and $k>n/2$. Then for each $R > 0$ there exist a $C > 0$ such
that
\eqn{Lip1}{ \norm{F(u_1)-F(u_2)}_{H^{k}_{\delta}} \leq
C\norm{u_1-u_2}_{H^{k}_{\delta}} \quad \text{for all $u_1,u_2\in
B_{R}(H^{k}_{\delta})$.} }
\end{lem}
\begin{proof} The proof of this lemma follows closely that
of Theorem 1, Section 5.5.2, in \cite{RunSic}. We will only prove
the case when $V \cong \Rbb$ and leave the general case to the
reader. Define
\eqn{Lip2}{ H(x,y) := \frac{F(x)-F(y)}{x-y} - F'(0) \, .}
Then $H(0,0)=0$ and $H\in C^{\ell}(\Rbb^2,\Rbb)$. Suppose
$u_1,u_2\in B_{R}(H^{k}_{\delta})$. Then by Lemma \ref{SobA}, there
exists a $C>0$ such that $\norm{u_\alpha}_{L^{\infty}_{\delta}} \leq
C\norm{u}_{H^{k}_{\delta}}$ for $\alpha=1,2$. Since $\delta \leq 0$,
we have $\norm{u_\alpha}_{L^{\infty}} \leq
C\norm{u_\alpha}_{L^{\infty}_{\delta}}$. Consequently,
\leqn{Lip3}{\norm{u_\alpha}_{L^{\infty}} \leq C_{1} \quad \alpha=1,2
}
for some constant $C_{1}>0$ independent of $u_1,u_2\in
B_{R}(H^{k}_{\delta})$. Let $\chi \in \Co(\Rbb)$ be a function such
that $0\leq \chi \leq 1$, $\chi|_{B_{C_{1}}(\Rbb)}=1$ and
$\text{supp}\, \chi \subset \subset B_{2C_{1}}(\Rbb)$. Define
$\Ht(x,y) := H(\chi(x),\chi(y))$. Then $\Ht(0,0)=0$ and $\Ht \in
C^{\ell}_{b}(\Rbb^{2},\Rbb)$. Also, \eqref{Lip3} implies that
$\Ht(u_1,u_2) = H(u_1,u_2)$ and hence
\eqn{Lip4}{ F(u_1)-F(u_2) = \Ht(u_1,u_2)(u_1-u_2) + F'(0)(u_1-u_2)
\, . }
Thus
\alin{Lip4}{ \norm{F(u_1)-F(u_2)}_{H^{k}_{\delta}} & \leq
C(\norm{\Ht(u_1,u_2)}_{H^{k}_{\delta}}+1)\norm{u_1-u_2}_{H^{k}_{\delta}}
&& \text{by Lemma
\ref{SobB}} \\
& \leq C\norm{u_1-u_2}_{H^{k}_{\delta}}  }
by Lemma \ref{MoserA} and equation \eqref{Lip3}. This completes the
proof.
\end{proof}

%################### moll.tex #############################33

\subsect{moll}{Mollifiers}

Let $j \in \Co(\Rbb^{n})$ be any function that satisfies $j\geq 0$,
$j(x)=0$ for $|x|\geq 1$,
and $\int_{\Rbb^{n}}j(x)\,d^n x =1$.
Following the standard prescription, we construct from $j$ the
mollifier $j_{\epsilon}(x) := \epsilon^{-n}j(x/\epsilon)$
 $(\epsilon > 0)$ and the smoothing operator
\eqn{smooth}{ J_{\epsilon}(u)(x) := j_{\epsilon}* u(x) =
\int_{\Rbb^{n}}j_{\epsilon}(x-y)u(y)\, d^{n}y \, .}
The next lemma shows that on the weighted spaces
$W^{k,p}_{\delta}$ the smoothing operator satisfies similar
estimates to the ones it satisfies on the unweighted spaces
$W^{k,p}$. This lemma is proved in the same manner as the Moser
estimates by using scaling and local estimates.
\begin{lem} \label{mollA} \mnote{mollA}
$\;$ \\ \vspace{-0.2cm}
\begin{itemize} \item[(i)] For $1\leq p \leq \infty$,
there exists a $C>0$
independent of $\epsilon \in (0,1)$ such that
\eqn{mollA1}{\norm{J_{\epsilon}u}_{W^{k,p}_{\delta}} \leq C
\norm{u}_{W^{k,p}_{\delta}}\, . }
\item[(ii)] For $1\leq p <
\infty$,
\eqn{mollA2}{\lim_{\epsilon\rightarrow 0^{+}}
\norm{J_{\epsilon}u-u}_{W^{k,p}_{\delta}} = 0 \, . }
\item[(iii)] For $1\leq p\leq \infty$,
$I\in \Nbb^{n}_{0}$, there exists a $C >0$
such that
\eqn{mollA3}{\norm{D^{I}J_{\epsilon}u}_{W^{k,p}_{\delta}} \leq
C\norm{u}_{W^{k,p}_{\delta}} \, . }
\item[(iv)] For $\ell > k$ there exists a $C>0$ independent
of $\epsilon \in (0,1)$ such that
\eqn{mollA3a}{ \norm{\Je u}_{H^k_\delta} \leq C
\epsilon^{\ell-k}\norm{u}_{H^{\ell}_\delta} \, .}
\item[(v)] For $1\leq p\leq \infty$ there exists a constant
$C>0$ independent of $\epsilon \in (0,1)$
such that
\eqn{mollA3b}{ \norm{\partial_{k}[f,\Je]u}_{L^p_{\delta_1+\delta_2}}
\leq C\norm{\nabla f}_{L^p_{\delta_1}}\norm{u}_{L^p_{\delta_2}} \, .
}
\end{itemize}
\end{lem}
\begin{proof}
\noindent \textbf{(i)-(iv)} On the
standard Sobolev spaces the smoothing operator
$\Je$ satisfies the well known properties:
\begin{itemize}
\item[(a)] For $1\leq p \leq \infty$,
there exists $C>0$
independent of $\epsilon> 0$ such that
\eqn{mollA4}{\norm{J_{\epsilon}u}_{W^{k,p}} \leq C
\norm{u}_{W^{k,p}} .}
\item[(b)] For $1\leq p < \infty$,
\eqn{mollA5}{\lim_{\epsilon\rightarrow 0^{+}}
\norm{J_{\epsilon}u-u}_{W^{k,p}} = 0 .}
\item[(c)] For $1\leq p\leq
\infty$, $I\in \Nbb^{n}_{0}$, there exists a $C >0$ such that
\eqn{mollA6}{\norm{D^{I}J_{\epsilon}u}_{W^{k,p}} \leq
C\norm{u}_{W^{k,p}} . }
\item[(d)] For $\ell > k$ there exists a constant $C>0$\
independent of $\epsilon \in (0,1)$ such that
\eqn{mollA7a}{ \norm{\Je u}_{H^{k}} \leq C
\epsilon^{\ell-k}\norm{u}_{H^{\ell}} .  }
\end{itemize}
%Proofs of (a)-(c) can be found in \cite{Adams} while
%(d) can be see most easily from the Fourier transform
%of $u$.
As in the previous section, the weighted estimates (i--iv) follow
from the standard ones (a--d) and scaling. We will only prove (i)
for $1\leq p < \infty$, and leave (ii--iv) to the reader.

{}From the definition of the partition of unity $\{\phi_j\}$, it
follows that
\lgath{mollA7}{ \phi_0(x)\sum_{k=0}^{\infty}\phi_k(x-y)
= \phi_0(x)\sum_{k=0}^{2}\phi_k(x-y)\, \label{mollA7.1} \\
\phi_1(x)\sum_{k=0}^{\infty}\phi_k(x-y)
=\phi_1(x)\sum_{k=0}^{3}\phi_k(x-y) \label{mollA7.2}\,
\intertext{and}
\phi_{j}(x)\sum_{k=0}^{\infty}\phi_{k}(x-y)
= \phi_j(x)\sum_{k=-2}^{k=+2} \phi_{j+k}(x-y) \quad j\geq 2
\label{mollA7.3}
}
for all $x\in \Rbb^n$ and $|y|\leq 1$. Fix $j\geq 2$ and suppose
$\epsilon < 1$. Then $j_{\epsilon}(x)=0$ for $|x|\leq \epsilon < 1$.
So \eqn{mollA8}{ \phi_j(x)\Je u(x) = \phi_j(x)\int_{\Rbb^n}
u(x-y)j_{\ep}(y)\,d^ny = \phi_j(x) \int_{|y|<1}u(x-y)j_{\ep}(y)
\,d^n y\, , } and thus
\alin{mollA9}{
\phi_j(x)\Je u(x)
&= \int_{|y|<1} \phi_j(x)\sum_{k=0}^{\infty}\phi_k(x-y)u(x-y)
j_{\ep}(y)\,d^n y  && \text{since ${\sum_j \phi_j =1}$}
\notag \\
& = \int_{|y|<1}\phi_j(x)\sum_{k=-2}^{2}\phi_{j+k}(x-y)
u(x-y)j_{\ep}(y)\,d^n y && \text{by \eqref{mollA7.3}} \notag \\
& = \phi_j(x)\Je\bigr(\sum_{k=-2}^{2}\phi_{j+k}u\bigl)(x)
\label{mollA9.1} \, .
}
So
\alin{mollA10}{
S_{j}(\phi_j\Je u) &= S_j\bigl(\phi_j\Je\bigl(\sum_{k=-2}^{2}\phi_{j+k}u
\bigr)\bigr) \\
& = \phi_1 J_{\ep/2^{j-1}}\bigl(\sum_{k=-2}^{2}S_{j}(\phi_{j+k}u)
\bigr) && \text{since $S_{j}\circ J_{\ep} = J_{\ep/2^{j-1}}\circ
S_j$}\\
& =  \phi_1 J_{\ep/2^{j-1}}\bigl(\sum_{k=-2}^{2}S_{-k+1}\circ
S_{j+k}(\phi_{j+k}u)
\bigr) \, .
}
Using this and inequality (a) yields
\leqn{mollA11}{ 2^{-p\delta(j-1)}\norm{S_j(\phi_j\Je u)}^p_{W^{k,p}}
\leq C \sum_{k=-2}^{2} 2^{-p\delta(j+k-1)}
\norm{S_{j+k}(\phi_{j+k}u)}^p_{W^{k,p}} }
for some constant $C>0$ independent of $\epsilon$. Similar
calculations show that
\leqn{mollA12}{ \norm{\phi_0\Je u}^p_{W^{k,p}}\leq
C\norm{u}^p_{W^{k,p}} \quad \text{and}\quad \norm{S_1(\phi_1\Je
u)}^p_{W^{k,p}}\leq C\norm{u}^p_{W^{k,p}} \, . }
It then follows from \eqref{scaleC}, \eqref{mollA11}, and
\eqref{mollA12} that $\norm{\Je u}_{W^{k,p}_{\delta}} \leq
C\norm{u}_{W^{k,p}_{\delta}}$ for some constant $C$ independent of
$\epsilon \in (0,1)$ which proves part (i).

\noindent\textbf{(v)} We will only prove the case $1<p<\infty$
and leave the rest to the reader. Let
\eqn{mollA13}{ v(x) := \int_{\Rbb^n}\partial_k
j_{\ep}(x-y)(f(x)-f(y))u(y)\, d^n y\, . } Then \leqn{mollA14}{
\partial_k\bigl([f,\Je]u)(x)
= \partial_{k}f(x)\Je(u)(x) + v(x) \, .}
Suppose $j\geq 2$. Then
\lalign{mollA15}{ \phi_j& (x)v(x) = \phi_j (x)
\int_{|y|<1}\frac{1}{\ep^{n}}\partial_k j\left(\frac{y}{\ep}\right)
\frac{f(x)-f(x-y)}{\ep}u(x-y)\, d^n y \notag \\
& = \int_{|y|<1}\frac{1}{\ep^{n}}
\partial_k j\left(\frac{y}{\ep}\right) \frac{f(x)-f(x-y)}{\ep}
\phi_j(x)\sum_{s=-2}^{2}\phi_{j+s}(x-y)u(x-y)\, d^n y \notag
&& \text{by \eqref{mollA7.2}} \\
&= \int_{\Rbb^n}\frac{1}{\ep^{n}}\partial_k j\left(\frac{x-y}{\ep}\right)
\frac{f(x)-f(x-y)}{\ep}\phi_j(x)\sum_{s=-2}^{2}\phi_{j+s}(y)u(y)\, d^n y\, .
\label{mollA15.1} \notag
}
Let
\eqn{mollA16}{ \psi_\ep(x,y) := \frac{1}{\ep^{n}}\partial_k
j\left(\frac{x-y}{\ep}\right) \frac{f(x)-f(x-y)}{\ep}\phi_j(x) \, .
} Then \eqn{mollA17}{ S_{j}(\phi_j v)(x) =
\int_{\Rbb^n}\psi_\ep(2^{j-1}x,y)
\sum_{s=-2}^{s=2}\phi_{j+s}(y)u(y)\, d^n y, }
and hence by H\"{o}lder's inequality
\leqn{mollA18}{ |S_{j}(\phi_j v)(x)| \leq \left(\int_{\Rbb^n}|
\psi_\ep(2^{j-1}x,y)|\, d^n y\right)^{1/p'} \left(\int_{\Rbb^n}
|\psi_\ep(2^{j-1}x,y)| \Bigl|
\sum_{s=-2}^{s=2}\phi_{j+s}(y)u(y)\Bigr|^p \, d^n y \right)^{1/p} }
where $1/p + 1/p' = 1$ .

Let
\eqn{mollA19}{ \At_j := \bigcup_{s=-2}^{0}A_{2^{j+s}}\, . }
Then
\eqn{mollA20}{ |f(x)-f(y)| \leq \norm{\nabla
f}_{L^\infty(\At_j)}|x-y| \quad \text{for all $x,y$ such that $x\in
A_{2^{j-1}}$ and $|x-y|\leq 1$.} }
Since $\partial_k j(x)=0$ for $|x|>1$ and $\text{supp}\, \phi_j
\subset A_{2^{j-1}}$, it follows that
\eqn{mollA21}{ |\psi_\ep(x,y)| \leq C\norm{\nabla
f}_{L^\infty(\At_j)} \omega_\ep(x-y)\, , } where \eqn{mollA22}{
\omega_\ep(x) := \frac{1}{\ep^n}|x/\epsilon|\,|\nabla
j(x/\epsilon)|\, . }
So
\eqn{mollA23}{ \int_{\Rbb^n}|\psi_\ep(x,y)|\, d^n y \leq
C\norm{\nabla f}_{L^\infty(\At_j)} \int_{\Rbb}\omega_{\ep}(x-y)\,
d^n y = C\norm{\omega_1}_{L^1}\norm{\nabla f}_{L^\infty(\At_j)} }
and the constant $C$ is independent of $\epsilon \in (0,1)$. Using
this in \eqref{mollA18} and integrating over $x$ yields
\lalign{mollA24}{ \norm{S_j(\phi_j v)}^p_{L^p} &\leq C \norm{\nabla
f}_{L^\infty(\At_j)}^{p/p'} \int_{\Rbb^n}\int_{\Rbb^n}
|\psi_\ep(2^{j-1}x,y)| \Bigl|
\sum_{s=-2}^{s=2}\phi_{j+s}(y)u(y)\Bigr|^p \, d^n y\, d^n x\notag \\
& \leq C\norm{\nabla f}_{L^\infty(\At_j)}^{p}
\int_{\Rbb^n}\int_{\Rbb^n}\omega_\ep(2^{j-1}x-y)
\Bigl|\sum_{s=-2}^{s=2}\phi_{j+s}(y)u(y)\Bigr|^p \, d^n y\, d^n x \notag \\
& = C\norm{\nabla f}_{L^\infty(\At_j)}^{p}
\int_{\Rbb^n}\Bigl|\sum_{s=-2}^{s=2}\phi_{j+s}(2^{j-1}y)u(2^{j-1}y)\Bigr|^p
\int_{\Rbb^n}\omega_\ep(2^{j-1}x-y)\, d^n x\, d^n y \notag \\
& \leq C\norm{\omega_1}_{L^1}\norm{\nabla f}_{L^\infty(\At_j)}^{p}
\int_{\Rbb^n}\Bigl|\sum_{s=-2}^{s=2}\phi_{j+s}(2^{j-1}y)u(2^{j-1}y)\Bigr|^p
\, d^n y \, . \notag
}
This shows that
\eqn{mollA25}{ 2^{-p(\delta_1+\delta_2)(j-1)} \norm{S_j(\phi_j
v)}^p_{L^p} \leq C 2^{-p\delta_1(j-1)}\norm{\nabla
f}_{L^\infty(\At_j)}^{p}
\sum_{s=-2}^{2}2^{-p\delta_2(j+s-1)}\norm{S_{j+s}(\phi_{j+s}u)}_{L^{p}}^p
}
for some constant $C>0$ independent of $\epsilon \in (0,1)$. Thus
\eqn{mollA26}{ 2^{-p(\delta_1+\delta_2)(j-1)} \norm{S_j(\phi_j
v)}^p_{L^p} \leq C \norm{\nabla f}_{L^\infty_{\delta_1}}
\sum_{s=-2}^{2}2^{-p\delta_2(j+s-1)}\norm{S_{j+s}(\phi_{j+s}u)}_{L^{p}}^p.
}
Similar arguments show that
\eqn{mollA26b}{ \norm{\phi_0 v}_{L^p} \leq C\norm{\nabla
f}_{L^\infty_{\delta_1}} \norm{u}_{L^{p}_{\delta_2}} \quad
\text{and}\quad \norm{S_1(\phi_1 v)}_{L^p} \leq C\norm{\nabla
f}_{L^\infty_{\delta_1}} \norm{u}_{L^{p}_{\delta_2}} }
where again the constant $C>0$ is independent of $\ep \in (0,1)$. It
then follows from \eqref{scaleC} that
\leqn{mollA27}{ \norm{v}_{L^{p}_{\delta_1+\delta_2}} \leq
C\norm{\nabla f}_{L^\infty_{\delta_1}} \norm{u}_{L^{p}_{\delta_2}}\,
. }
Applying this to \eqref{mollA14} yields
\alin{mollA27}{
\norm{\partial_k[f,\Je]u}_{L^{p}_{\delta_1+\delta_2}}
&\leq \norm{\nabla f}_{L^{p}_{\delta_1}}\norm{\Je u}_{L^p_{\delta_2}}
+  C\norm{\nabla f}_{L^\infty_{\delta_1}}
\norm{u}_{L^{p}_{\delta_2}} \\
& \leq C\norm{\nabla f}_{L^\infty_{\delta_1}}
 \norm{u}_{L^{p}_{\delta_2}} && \text{by part (i)}
}
for some constant $C$ independent of $\epsilon \in (0,1)$.
\end{proof}

The operator $\Je$ is no-longer self-adjoint on the weighted
spaces $L^2_\delta$. However, the next lemma shows
that the difference between $\Je$ and its adjoint is
a operator that smooths one derivative and its operator
norm can be estimated independent of $\epsilon$. This
will turn out to be crucial in proving the existence
and uniqueness of solutions to quasilinear parabolic
equations using our Galerkin method.

\begin{lem} \label{mollB} \mnote{mollB}
Let $\Jde$ be the adjoint of $\Je$ on $L^2_\delta$ with respect to
the inner product \eqref{L2ip} and $\rho = \sigma^{-2\delta-n}$.
Then $\Jde = \Je + \rho^{-1}[\Je,\rho]$,
\eqn{mollB1}{ \norm{\Jde u}_{L^2_\delta} \leq C\norm{u}_{L^2_\delta}
\, , } and \eqn{mollB2}{ \ip{\partial_k\Jde u}{v}_{L^2_\delta} \leq
\ip{\partial_k\Je u}{v}_{L^2_\delta} +
C\norm{u}_{L^2_\delta}\norm{v}_{L^2_\delta} }
for some constant $C>0$ independent of $\epsilon\in (0,1)$.
\end{lem}
\begin{proof}
First we note that from standard properties of adjoints, we have
$\norm{\Je^{\dagger}}_{op}=\norm{\Je}_{op}$. Lemma \ref{mollA}.(i)
shows that $\norm{\Je}_{op}\leq C$ for some $C$ independent of
$\epsilon$. Thus $\norm{\Je^{\dagger}}_{op}\leq C$ and hence
$\norm{\Je^{\dagger}u}_{L^{2}_{\delta}}\leq
C\norm{u}_{L^{2}_{\delta}}$ for all $u \in L^{2}_{\delta}$. From the
definition of the adjoint, we have $\ip{\Jde u}{v}_{L^2_\delta} =
\ip{u}{\Je v}_{L^2_\delta}$ for all $u$, $v$ $\in L^2_\delta$. So
\alin{mollB4}{
\ip{\Jde u}{v} &= \ip{\rho u}{\Je v}_{L^2}
= \ip{\Je(\rho u)}{v}_{L^2}
= \ip{\rho^{-1}\Je(\rho u)}{v}_{L^2_\delta}
= \ip{(\Je + \rho^{-1}[\Je,\rho])u}{v} \, .
}
Since $u,v \in L^2_\delta$ where arbitrary, this proves that $\Jde =
\Je + \rho^{-1}[\Je,\rho]$. Therefore using Lemma \ref{mollA} and
the Cauchy-Schwartz and weighted H\"{o}lder inequalities, we get
\alin{mollB5}{
\ip{\partial_k&\Jde u}{v}_{L^2_\delta}
= \ip{\partial_k\Je u}{v}_{L^2_\delta}
-\ip{\rho^{-2}\partial_k\rho\Je(\rho u)}{v}_{L^2_\delta}
+\ip{\rho^{-1}\partial_k\rho\Je u}{v}_{L^2_\delta}
+\ip{\rho^{-1}\partial_k[\Je,\rho]u}{v}_{L^2_\delta} \\
& \leq  \ip{\partial_k\Je u}{v}_{L^2_\delta}
+ \bigl(\norm{\rho^{-2}\partial_k\rho}_{L^\infty_{\delta
+2\delta+n}}\norm{\Je(\rho u)}_{L^2_{\delta-2\delta-n}}
+ \norm{\rho^{-1}\partial_k\rho}_{L^\infty}
\norm{\Je u}_{L^2_{\delta}} \\
& \qquad \qquad +
 \norm{\rho^{-1}}_{L^\infty_{2\delta+n}}
\norm{\partial_k[\Je,\rho]u}_{L^2_{\delta-2\delta-n}} \bigr)
\norm{v}_{L^2_\delta}
 \\
&\leq C\bigl(\norm{\rho u}_{L^2_{\delta-2\delta-n}}
+ (1+\norm{\rho}_{L^\infty_{-2\delta-n}})
\norm{u}_{L^2_{\delta}}\bigr)\norm{v}_{L^2_\delta}
 \leq C\norm{u}_{L^2_\delta}\norm{v}_{L^2_\delta}
}
where the constant $C$ is independent of $\ep \in (0,1)$.
\end{proof}

%################### par.tex ####################################3

\sect{par}{Parabolic equations on the $H^{k}_{\delta}$ spaces}
\setcounter{equation}{0}

In this section, we will prove a local existence theorem for
parabolic equations on the $H^{k}_{\delta}$ spaces. We do this by
adapting the local existence proof of Taylor \cite{TayIII} (see
Theorem 7.2, pg 330, and Proposition 7.7, pg 334) which is valid for
the $H^{k}$ spaces. The parabolic equations that we will consider
are of the form
\lalign{parabolic}{ &\partial_{t}u = a^{ij}(v,u)\partial_i\partial_j
u + b(v,u,\nabla u)
+ f, \label{parabolic.1}\\
&u|_{t=0}  = u_{0} \label{parabolic.2}}
where
\begin{itemize}
\item[(i)] the maps $u=u(t,x)$ and $f=f(t,x)$ are $\Rbb^{m}$-valued
while $v=v(t,x)$ is a $\Rbb^{r}$-valued,
\item[(ii)] $a^{ij}\in
C^{\infty}(\Rbb^{m+r}, \mathbb{M}_{m\times m})$,
$a^{ij}$ is symmetric for each $i,j=1,\ldots,n$,
\item[(iii)] $b\in C^{\infty}(\Rbb^{r+m(1+n)},\mathbb{M}_{m\times m})$, $b(0)=0$, and
\item[(iv)] there exists a constant $\omega > 0$ such that
\leqn{uellipt}{ \ipe{w}{a^{ij}(u,v)\xi_i\xi_j\cdot w} \geq \omega
\enorm{\xi}^2 \enorm{w}^2 \quad \text{for all $u,w\in \Rbb^{m}$,
$v\in \Rbb^r$ and  $\xi \in \Rbb^{n}$.} }
\end{itemize}

\subsect{galerk}{Galerkin method}

Following Taylor (\cite{TayIII}, Sec.~15.7), we first solve the
approximating equation
\lalign{approx}{ & \partial_t\ue =  \Je a^{ij}(v,\Je
\ue)\partial_i\partial_j \Je \ue +  \Je b(v,\Je \ue,\nabla \Je
\ue) + \Je f \label{approx.1} \\
& \ue|_{t=0} = u_{0}, \label{approx.2} }
and latter show that the solutions $u_\ep$ converge to a solution of
\eqref{parabolic.1}--\eqref{parabolic.2} as $\ep \rightarrow 0$.
\begin{prop} \label{parA} \mnote{[parA]} Suppose $\epsilon > 0$,
$\delta \leq 0$, $\ell\geq k>n/2$, $u_{0}\in H^{k}_{\delta}$, and
$v,f\in C^{0}([0,T],H^{k}_{\delta})$ for some $T>0$. Then there
exists a $T_{*}\in (0,T]$ and a unique $u_{\epsilon}\in
C^{1}([0,T_{*}),H^{k}_{\delta})$ that solves the initial value
problem \eqref{approx.1}--\eqref{approx.2}. Moreover if $\sup_{0\leq
t < T_{*}(\epsilon)}$ $\norm{u_{\epsilon}(t)}_{H^k_\delta}$ $<
\infty$ then there exists a $T^{*} \in (T_{*},T)$ such that $\ue$
extends to a unique solution on $[0,T^{*}]$.\end{prop}
\begin{proof} Let $R=\max\{\norm{u_{0}}_{H^{k}_{\delta}},
\sup_{ 0\leq t\leq T}\norm{v(t)}_{H^{k}_{\delta}}\}+1 $ and
\eqn{parA1}{F(t,w) := \Je a^{ij}(v,\Jde w)\partial_i\partial_j \Je w
+ \Je b(v,\Je w,\nabla \Je w) + \Je f \, .}
Then the approximating equations \eqref{approx.1}--\eqref{approx.2}
can be written as the first order differential equation $\dot{u} =
F(t,u) \; ; \; u(0)= u_{0}$ on $H^{k}_{\delta}$. It follows from
Lemmata \ref{Lip}, \ref{mollA}, and \ref{mollB} that $F\in
C^{0}([0,T]\times B_{R}(H^{k}_{\delta}),H^{k}_{\delta} )$ and also
that there exists a constant $C > 0$ such that
$\norm{F(t,w_{1})-F(t,w_{2})}_{H^{k}_{\delta}} \leq
C\norm{w_1-w_2}_{H^{k}_{\delta}}$ for all $w_1,w_2\in
B_{R}(H^{k}_{\delta})$. Therefore we can apply the standard
existence theorem for first order equations on Banach spaces (see
\cite{AMR}, Lemma 4.1.6) to conclude that there exists a $T_{*} \in
(0,T]$ such that the initial value problem $\dot{u} = F(t,u)\; ; \;
u(0)= u_{0}$ has a unique solution $u \in
C^{1}([0,T_{*}),H^{k}_{\delta})$. Also, standard ODE results show
that the solution can be extended as long as $u$ is bounded.
\end{proof}

\subsect{energy}{Energy estimates}

We will now assume that $k>n/2+1$. By Proposition \ref{parA}, we
have a sequence of solutions $u_{\epsilon} \in
C^{1}([0,T(\epsilon)],H^{k}_{\delta})$ to the approximating equation
\eqref{approx.1}--\eqref{approx.2}. Setting $\ute:= \Je\ue$ and
differentiating $\norm{\ue}^2_{H^{k}_\delta}$ with respect to $t$
yields
\leqn{eng1}{ \frac{d\:}{dt} \frac{1}{2}\norm{\ue}^2_{H^k_\delta} =
\ip{\ue}{\Je a^{ij}_{\ep}\partial_i\partial_j \ute}_{H^k_\delta} +
\ip{\ue}{\Je b_\ep}_{H^k_\delta} +\ip{\ue}{\Je f}_{H^k_\delta}\, , }
where $a^{ij}_\ep := a^{ij}(v,\ute)$ and $b_\ep := b(v,\ute,\nabla
\ute)$. Since $\delta \leq 0$ and $k>n/2+1$, the weighted Sobolev
inequality (Lemma \eqref{SobA}) implies that
\leqn{eng2a}{\norm{w}_{W^{1,\infty}} \leq
\norm{w}_{W^{1,\infty}_{\delta}} \leq C \norm{w}_{H^{k}_{\delta}}\,
. }

The Cauchy-Schwartz inequality and Lemma \ref{MoserA} show that
\leqn{eng2}{ \ip{\ue}{\Je f}_{H^k_\delta} \leq
C\norm{\ue}_{H^{k}_{\delta}}} and \leqn{eng3}{ \ip{ \ue}{\Je
b_\ep}_{H^k_\delta} \leq p(\norm{
\ue}_{W^{1,\infty}})(1+\norm{\ue}_{H^k_\delta}+\norm{\nabla
\ute}_{H^k_\delta}) \norm{\ue}_{H^k_\delta}\,  }
where we use the notation $p(x)$ to denote a polynomial that is
independent of $\epsilon$.

Let $\rho = \sigma^{-2(\delta-|I|)-n}$. Then
\alin{eng4}{
\ip{D^I \ue}{\Je a^{ij}_\ep & D^{I}\partial_{i}\partial_j
\ute}_{L^2_{\delta-|I|}}  = -\ip{\partial_i \Jde D^I u}{
a^{ij}_\ep D^I\partial_j \ute}_{L^2_{\delta-|I|}}    - \notag \\ & \ip{(\partial_{i} a^{ij}_\ep +
\rho^{-1}\partial_i \rho)\Jde D^I u}{D^I\partial_j\ute}_{L^2_{\delta-|I|}}
\leq -\omega \norm{D^I \nabla \ute}^2_{L^2_{\delta-|I|}} + \notag \\ &
(\norm{a_{\ep}}_{W^{1,\infty}} +
\norm{\rho^{-1}\nabla \rho}_{L^{\infty}})\norm{D^I
u}_{L^2_{\delta-|I|}}\norm{D^{I}\nabla \ute}_{L^2_{\delta-|I|} }
}
where in deriving the last inequality we used \eqref{uellipt} and
Lemma \ref{mollB}. Lemma \ref{mollA} implies that
$\norm{a_\ep}_{W^{1,\infty}}\leq C(1+\norm{\ue}_{W^{1,\infty}})$.
Using this and $\norm{\rho^{-1}\nabla \rho}_{L^{\infty}} <\infty$,
the above inequality implies that
\lalign{eng5aa}{ \ip{D^I\ue}{\Je a^{ij}_\ep
D^I\partial_{i}\partial_j \ute}_{L^2_{\delta-|I|}} \leq & -\omega
\norm{D^I \nabla
\ute}^2_{L^2_{\delta-|I|}} + \notag \\
& p(\norm{\ue}_{W^{1,\infty}})\norm{D^I\ue}_{L^2_{\delta-|I|}}
\norm{D^{I}\nabla\ute}_{L^2_{\delta-|I|}} \, . \label{eng5} }
Also,
\alin{eng6}{ \ip{D^I&\ue}{\Je
[D^I,a^{ij}_{\ep}]\partial_i\partial_j\ute}_{L^2_{\delta-|I|}} \leq
\norm{D^I\ue}_{L^2_{\delta-|I|}}\norm{[D^I,\Je a^{ij}_{\ep}]
\partial_i\partial_j\ute}_{L^2_{\delta-|I|}} \\
& \leq C \norm{\ue}_{H^{k}_{\delta}}\bigl(\norm{\nabla
a_{\ep}}_{H^{k-1}_{\delta-1}}
\norm{\partial_{i}\partial_j\ute}_{L^{\infty}}+
\norm{Da_\ep}_{L^{\infty}}\norm{\partial_{i}\partial_j\ute}_{H^{k-1}_{\delta-1}}\bigr)
}
by Lemmata \ref{MoserC} and \ref{mollA}. From the definition of the
weighted norm, we get that
$\norm{\partial_{i}\partial_j\ute}_{H^{k-1}_{\delta-1}}\leq
C\norm{\nabla \ute}_{H^{k}_{\delta}}$ while \eqref{eng2a} shows that
$\norm{\partial_i\partial_j \ute}_{L^{\infty}}\leq
C\norm{\nabla\ute}_{H^{k}_{\delta}}$. Moreover, $\norm{\nabla
a_{\ep}}_{H^{k-1}_{\delta-1}} \leq
p(\norm{\ue}_{L^{\infty}})(1+\norm{\ue}_{H^{k}_{\delta}})$ by
Lemmata \ref{MoserA} and \ref{mollA}. Therefore
\leqn{eng7}{ \ip{D^I\ue}{\Je
[D^I,a^{ij}_{\ep}]\partial_i\partial_j\ute}_{L^2_{\delta-|I|}}
\leq p(\norm{\ue}_{W^{1,\infty}})(1+\norm{\ue}_{H^k_\delta}+
\norm{\nabla\ute}_{H^k_\delta})\norm{\ue}_{H^k_\delta} \,. }
Adding the two inequalities \eqref{eng5} and \eqref{eng7} and then
summing over $0\leq |I|\leq k$ yields
\leqn{eng8}{\ip{\ue}{a^{ij}_{\ep}\partial_i\partial_j\ute}_{H^k_\delta}
\leq -\omega\norm{\nabla\ute}^2_{H^k_\delta} +
p(\norm{\ue}_{W^{1,\infty}})(1+\norm{\ue}_{H^k_\delta}+
\norm{\nabla \ute}_{H^k_\delta})\norm{\ue}_{H^k_\delta}\, . }
{}From \eqref{eng1}-\eqref{eng3} and \eqref{eng8}, we get
\leqn{eng9}{ \frac{d\:}{dt}\frac{1}{2}\norm{\ue}^2_{H^k_\delta} \leq
 -\omega\norm{\nabla \ute}^2_{H^k_\delta} +
p(\norm{\ue}_{W^{1,\infty}})(1+\norm{\ue}_{H^k_\delta}+
\norm{\nabla\ute}_{H^k_\delta})\norm{\ue}_{H^k_\delta}\, . }
Using $AB\leq \omega A^2 + (1/4\omega)B^2$ with $A=\norm{\nabla
\ute}_{H^k_\delta}$ and
$B=p(\norm{\ue}_{W^{1,\infty}})\norm{\ue}_{H^k_\delta}$, yields
\leqn{eng10}{ \frac{d\:}{dt}\frac{1}{2}\norm{\ue}^2_{H^k_\delta}
\leq p(\norm{\ue}_{W^{1,\infty}})(1+\norm{\ue}_{H^k_\delta})
\norm{\ue}_{H^k_\delta}\, .}
Finally, using \eqref{eng2a}, we arrive at
\leqn{engll}{ \frac{d\:}{dt}\norm{\ue}_{H^k_\delta} \leq
p(\norm{\ue}_{H^k_\delta})(1+\norm{\ue}_{H^k_\delta})\, . }
Then Gronwall's inequality implies that there exists a constant
$C>0$ and a $T_{*}\in (0,T)$, both independent of $\ep$, such that
$T(\epsilon)\geq T_{*}$ for all $\ep > 0$, and
\leqn{eng12}{ \norm{\ue(t)}_{H^{k}_\delta} \leq C \quad \text{for
all $t \in [0,T_{*})\,$.} }
Also, \eqref{approx.1} and Lemma \ref{mollA} imply that
\eqn{eng13}{\norm{\partial_t \ue}_{H^{k-2}_\delta} \leq
C\norm{a_\ep^{ij}\partial_{i}\partial_{j}\ute}_{H^{k-2}_\delta} +
C\norm{b_\ep}_{H^{k-2}_\delta} + C\norm{f}_{H^{k-2}_\delta} }
and hence
\leqn{eng14}{ \norm{\partial_t \ue(t)}_{H^{k-2}_\delta} \leq C \quad
\text{for all $t\in [0,T_{*})$} }
by \eqref{eng12} and Lemmata \ref{MoserA} and \ref{mollA} where
again the constant $C$ is independent of $\ep$.

\subsect{loc}{Local existence}

We are now ready to prove local existence of solutions to
\eqref{parabolic.1}--\eqref{parabolic.2}. The following theorem is
the weighted version of Theorem 7.2, pg 330 in \cite{TayIII}.
\begin{thm} \label{locA} \mnote{[locA]} Suppose$\delta \leq 0$,
$\ell\geq k>n/2+1$, $u_{0}\in H^{k}_{\delta}$, and $v,f\in
C^{0}([0,T],H^{k}_{\delta})$ for some $T>0$. Then there exists a
$T_{*}\in (0,T)$ and a $u \in L^{\infty}((0,T_{*}),H^{k}_{\delta})
\cap \emph{\text{Lip}}([0,T_{*}),H^{k-2}_{\delta})$ that solves the
initial value problem \eqref{parabolic.1}--\eqref{parabolic.2}.
\end{thm}
\begin{proof} In the previous section we established that
$\ue \in C^0([0,T^{*}],H^k_\delta)\cap
C^1([0,T^{*}],H^{k-2}_\delta)$ is uniformly bounded for some
$T_{*}\in (0,T)$. But, $C^0([0,T_{*}],H^k_\delta)\cap
C^1([0,T_{*}],H^{k-2}_\delta) \subset
L^\infty((0,T_{*}),H^k_\delta)\cap
W^{1,\infty}((0,T_{*}),H^{k-2}_\delta)$ and
$L^\infty((0,T_{*}),H^k_\delta)$ and
$W^{1,\infty}((0,T_{*}),H^{k-2}_\delta)$ are the dual of a Banach
space which implies via the Banach-Alaoglu theorem that from any
bounded sequence we can extract a subsequence that converges in the
weak* topology.
\footnote{$W^{1,\infty}((0,T_{*}),H^{k-2}_{\delta}) =
\text{Lip}([0,T_{*}),H^{k-2}_{\delta})$ by Theorem 2, pg 286, in
\cite{Ev}.}
Therefore there exists a sequence $\{\ep_{n}\}\subset (0,1)$ with
$\lim_{n\rightarrow \infty}\ep_{n}=0$ such that $u_n :=u_{\ep_n}
\rightarrow u\in L^{\infty}((0,T_{*}),H^k_\delta)\cap
W^{1,\infty}((0,T_{*}),H^{k-2}_{\delta})$ in the weak* topology as
$n\rightarrow \infty$.

The Sobolev interpolation inequality implies that
\eqn{interp}{\norm{v}_{H^{k-s}(B_{2^m})} \leq
C\norm{v}_{L^2(B_{2^m})}^{1-(k-s)/k}\norm{v}_{H^{k}(B_{2^m})}^{(k-s)/k}
}
for $s \in [0,k]$ and $m\in \Nbb_{0}$. Since
$H^{k-s}_{\delta}(B_R)\cong H^{k-s}(B_R)$, it follows that for
$s>0$, the sequence $u_n$ is bounded in
$C^{\sigma}([0,T_{*}],H^{k-s}(B^{2^m}) ) $ for some $\sigma > 0$.
The compactness of the imbedding of  $H^{k-s}(B_R)$ in
$H^{k-2-s_1}(B_R)$ $(s_1>s)$ then implies via the Ascoli theorem
that there exist a subsequence $u_{n_{m}}$ converging strongly to
$u$ in $C^{0}([0,T_{*}],H^{k-s}(B_{2^m}))$ for any $s>0$. Taking the
diagonal subsequence, we see that there exists a subsequence of
$u_n$, which we will again denote by $u_n$, that converges strongly
to $u$ in $C^{0}([0,T_{*}],H^{k-s}(B_{2^m}))$ for all $m\in \Nbb$.
For $s$ small enough,  $k-s>n/2+1$ and hence  the standard Moser
estimates imply that the map $H^{k-s}(B_{2^m})\times
H^{k-s}(B_{2^m}) \ni (w_1,w_2) \mapsto a^{ij}(w_1,w_2)\in
H^{k-s}(B_{2^m})$ is continuous. So finally the map
$H^{k-s}(B_{2^m})\times H^{k-s}(B_{2^m}) \times H^{k-s}(B_{2^m})\ni
(w_1,w_2,w_3) \mapsto a^{ij}(w_1,w_2)\partial_i\partial_j w_3 \in
H^{k-s-2}(B_{2^m})$ is continuous by the multiplication lemma. The
same arguments show that the map $H^{k-s}_{\delta}(B_{2^m})\times
H^{k-s}_{\delta}(B_{2^m}) \ni (w_1,w_2) \mapsto b(w_1,w_2,\nabla
w_2)\in H^{k-s-2}_{\delta}(B_{2^m})$ is also continuous. Therefore
\eqn{locA3}{ \Je a^{ij}_{\ep_n}
\partial_{i}\partial_j u_n + \Je b_{\ep_n} + \Je f
\longrightarrow a^{ij}(v,u)
\partial_{i}\partial_j u + b(v,u,\nabla u) + f }
strongly in $L^2(B_{2^m})$ as $n\rightarrow \infty$ for each $m\in
\Nbb_0$. This and the fact that $\partial_t u_n \rightarrow
\partial_t u$ weak* implies that
$\partial_t u =  a^{ij}(v,u)
\partial_{i}\partial_j u + b(v,u,\nabla u) + f$ a.e. on $(0,T_*)\times \Rbb^n$.
\end{proof}
Using the estimates of Sections \ref{Moser}, \ref{moll}, and
\ref{energy}, it is not difficult to adapt the proofs  of
Proposition 7.3--7.7, pp.~332--334 in \cite{TayIII} to get the
following theorem.
\begin{thm} \label{locB} \mnote{[locB]} The solution $u\in
L^{\infty}([0,T_{*}],H^k_\delta)\cap\emph{\text{Lip}}([0,T_{*}],
H^{k-2}_{\delta})$ from Theorem \ref{locA} is unique and satisfies
the additional regularity
\eqn{locB1}{u \in C^{0}([0,T_{*}),H^k_\delta)\cap
C^{1}([0,T_{*}),H^{k-2}_{\delta}) \cap
C^{0}([T_1,T_2],H^\ell_\delta) \cap
C^{1}([T_1,T_2],H^{\ell-2}_\delta)\, }
for every closed interval $[T_1,T_2] \subset (0,T_{*})$. Moreover,
if $\;\sup_{0\leq t < T_{*}}\norm{u(t)}_{W^{1,\infty}} < \infty$
then there exist a $T^{*} \in (T_{*},T)$ such $u$ can be extended to
a solution of \eqref{parabolic.1}--\eqref{parabolic.2} on
$[0,T^{*})$.\end{thm}

With more information on the structure of the function $b$,
it is possible to relax the requirement that
$v\in C^{0}([0,T),H^k_\delta)$ to $v\in C^{0}([0,T),H^k_\eta)$
for any $\eta \leq 0$ independent of $\delta$. Essentially
what this requires is that $b$ is of the form
$b(v,u,\nabla u) = b^{i}(v,u,\nabla u)\partial_i u
+ c(v,u,\nabla u) u$. A simple example of
this is dealt with in the next theorem where we consider
linear equations. The proof, which we omit, requires
only small changes to the above arguments.

\begin{thm} \label{locC} \mnote{locC}
Suppose, $\delta,\eta \leq 0$, $\ell\geq k>n/2+1$,
$u_{0}\in H^{k}_{\delta}$, and $v,f^i,c\in
C^{0}([0,T),H^{\ell}_{\eta})$ for some $T>0$. Then the initial
value problem
\eqn{locC1}{
\partial_t v = a^{ij}(v)\partial_i\partial_j v
+ b^i\partial_i u + c u = 0 \quad : \quad u(0)=u_0
}
has a unique solution
$u \in C^{0}([0,T),H^{k}_{\delta})
\cap C^{1}([0,T),H^{k-2}_{\delta})$.
Moreover, for each interval $[T_1,T_2]\subset (0,T)$, $u$
satisfies
$u \in C^{0}([T_1,T_2],H^{\ell}_\delta)
\cap C^{1}([T_1,T_2],H^{\ell-2}_{\delta})$.
\end{thm}

\end{document}